\documentclass[11pt]{article}
\usepackage{amsfonts}
\usepackage{amsmath}

\newtheorem{theorem}{Theorem}

\newtheorem{corollary}[theorem]{Corollary}

\newtheorem{definition}[theorem]{Definition}

\newtheorem{lemma}[theorem]{Lemma}

\newtheorem{proposition}[theorem]{Proposition}
\newtheorem{remark}[theorem]{Remark}

\newenvironment{proof}[1][Proof]{\noindent\textbf{#1.} }{\ \rule{0.5em}{0.5em}}
\input{tcilatex}

\begin{document}

\title{A topological theory of Maslov indices for Lagrangian and symplectic
paths}
\author{M. de Gosson and S. de Gosson \\
maurice.degosson@gmail.com}
\date{}
\maketitle

\begin{abstract}
We propose a topological theory of the Maslov index for lagrangian and
symplectic paths based on a minimal system of axioms. We recover, as
particular cases, the H\"{o}rmander and the Robbin--Salomon indices.
\end{abstract}

\section{Introduction}

The theory of the Maslov index has recently found new applications in
intersection theory.

\subsection{General notations\label{not}}

We essentially use the notations of Leray \cite{Leray}. Let $X=\mathbb{R}%
^{n} $; the product $X\times X^{\ast }$ is endowed with the canonical
symplectic structure defined by: 
\begin{equation*}
\omega (z,z^{\prime })=\left\langle p,x^{\prime }\right\rangle -\left\langle
p^{\prime },x\right\rangle
\end{equation*}%
if $z=(x,p)$, $z^{\prime }=(x^{\prime }$,$p^{\prime })$. The symplectic
group of $(X\times X^{\ast },\omega )$ will be denoted by $Sp(n)$: we have $%
s\in Sp(n)$ if and only if $\omega (sz,sz^{\prime })=\omega (z,z^{\prime })$
for all $z,z^{\prime }$. The Poincar\'{e} group $\pi _{1}(Sp(n))$ is
isomorphic to $(\mathbb{Z},+)$ and $Sp(n)$ has thus covering groups $%
Sp_{q}(n)$ of all orders $q=1,2,...,\infty $ ($Sp_{\infty }(n)$ is the
universal covering group of $Sp(n)$). The unitary group $U(n)$ is identified
with a subgroup of $Sp(n)$ by defining the action of $u=A+iB$ on $z=x+p$ by
the formula 
\begin{equation*}
\left( 
\begin{array}{c}
x^{\prime } \\ 
p^{\prime }%
\end{array}%
\right) =\left( 
\begin{array}{cc}
A & -B \\ 
B & A%
\end{array}%
\right) \left( 
\begin{array}{c}
x \\ 
p%
\end{array}%
\right)
\end{equation*}

We denote by $\Lambda (n)$ the Lagrangian Grassmannian of $(X\times X^{\ast
},\omega )$: $\ell \in \Lambda (n)$ if and only if $\ell $ is an $n$
dimensional subspace of $X\times X^{\ast }$ on which the form $\omega $
vanishes. The manifold $\Lambda (n)$ is connected and compact manifold and
has dimension $n(n+1)/2$. Its Poincar\'{e} group $\pi _{1}(\Lambda (n))$ is
isomorphic to $(\mathbb{Z},+)$ and $\Lambda (n)$ thus has coverings $\Lambda
_{q}(n)$ of all orders $q$. Its universal covering $\Lambda _{\infty }(n)$
is sometimes called the \textquotedblleft Maslov bundle\textquotedblright .

The group $Sp(n)$, and hence $U(n)$, acts on $\Lambda (n)$. The action of $%
U(n)$, and hence that of $Sp(n)$, is transitive: the orbit of every $\ell
\in \Lambda (n)$ under the action of $Sp(n)$ is $\Lambda (n)$.

We shall denote by $St_{\ell }(n)$ the stabilizer group of a Lagrangian $%
\ell $ in $Sp(n)$: 
\begin{equation*}
St_{\ell }(n)=\left\{ s\in Sp(n):s\ell =\ell \right\} \text{.}
\end{equation*}%
One can show (see \cite{Arnold,Dazord,Leray}) that the fibration $%
Sp(n)/St_{X^{\ast }}(n)=\Lambda (n)$ defines an isomorphism 
\begin{equation*}
\mathbb{Z}\cong \pi _{1}\left[ Sp(n)\right] \longrightarrow \pi _{1}\left[
\Lambda (n)\right] \cong \mathbb{Z}
\end{equation*}%
which is multiplication by $2$ on $\mathbb{Z}$. It follows (see \cite{Leray}%
, Theorem 3,3$^{o}$), p.36) that the action of $Sp(n)$ on $\Lambda (n)$ can
be lifted to a transitive action of the universal covering $Sp_{\infty }(n)$
on the Maslov bundle $\Lambda _{\infty }(n)$ such that 
\begin{equation}
(\alpha s_{\infty })\ell _{\infty }=\beta ^{2}(s_{\infty }\ell _{\infty
})=s_{\infty }(\beta ^{2}\ell _{\infty })  \label{alfabeta}
\end{equation}%
for all $(s_{\infty },\ell _{\infty })\in Sp_{\infty }(n)\times \Lambda
_{\infty }(n)$ ; $\alpha $ (resp. $\beta $) is the generator of $\pi
_{1}(Sp(n))$ (resp. $\pi _{1}(\Lambda (n))$) whose image in $\mathbb{Z}$ is $%
+1$.

The Lagrangian Grassmannian $\Lambda (n)$ is a stratified manifold; for $%
\ell \in \Lambda $ and $k$ an integer, $0\leq k\leq n$, 
\begin{equation*}
\Lambda _{\ell }(n,k)=\left\{ \ell ^{\prime }\in \Lambda (n):\dim (\ell \cap
\ell ^{\prime })=k\right\}
\end{equation*}%
is the stratum of $\Lambda (n)$ of order $k$, relative to the Lagrangian
plane $\ell $. We will also use the following notations: 
\begin{equation*}
\Lambda ^{2}(n;k)=\bigcup_{\ell \in \Lambda }\Lambda _{\ell }(k)\text{ \ , \ 
}\Lambda _{\infty }^{2}(n;k)=\pi ^{-1}(\Lambda ^{2}(k))
\end{equation*}%
that is: 
\begin{equation*}
\Lambda ^{2}(n;k)=\left\{ (\ell ,\ell ^{\prime })\in \Lambda (n)^{2}:\dim
(\ell \cap \ell ^{\prime })=k\right\}
\end{equation*}%
\begin{equation*}
\Lambda _{\infty }^{2}(n;k)=\left\{ (\ell _{\infty },\ell _{\infty }^{\prime
})\in \Lambda _{\infty }(n)^{2}:\dim (\ell \cap \ell ^{\prime })=k\right\} 
\text{.}
\end{equation*}%
In particular, $\Lambda _{\infty }^{2}(n;0)$ is the set of pairs $(\ell
_{\infty },\ell _{\infty }^{\prime })$ of elements of $\Lambda _{\infty }(n)$
whose projections $\ell $ and $\ell ^{\prime }$ are transversal. One proves
(see for instance \cite{Treves}) that $\Lambda _{\ell }(n;0)$ is an open
subset of $\Lambda (n)$ and that the $\Lambda _{\ell }(n;k)$ are, for $0\leq
k\leq n$, connected submanifolds of $\Lambda (n)$, of codimension $k(k+1)/2$%
. The closed set 
\begin{equation*}
\Sigma _{\ell }=\Lambda (n)\setminus \Lambda _{\ell
}(n;0)=\bigcup_{k=1}^{n}\Lambda _{\ell }(k)
\end{equation*}%
is called the \textquotedblleft Maslov cycle relative to $\ell $%
\textquotedblright . It is the set of Lagrangians that are not transverse to 
$\ell $. When $\ell =X^{\ast }$ we call $\Sigma _{\ell }$ the \textit{Maslov
cycle}, and denote it by $\Sigma $.

\subsection{Cohomological notations}

It will be economical to use the following elementary notations from the
singular Alexander--Spanier cohomology: let $E$ be a set, $k$ a positive
integer or $0$, and $(G,+)$ an abelian group. We shall call the functions $%
f:E^{k+1}\longrightarrow G$ $k$-\textit{cochains on} $E$ \textit{with values
in} $G$ (or simply $k$-\textit{cochains}$,$ or \textit{cochains}). The
coboundary $\partial f$ of a $k$-cochain is the $(k+1)$-cochain defined by: 
\begin{equation*}
\partial f(a_{0},...,a_{k+1})=\sum_{j=0}^{k+1}(-1)^{j}f(a_{0},...,\widehat{a}%
_{j},...,a_{k+1})\text{,}
\end{equation*}%
where the cap \symbol{94}\ deletes the term it covers. We evidently have
that $\partial ^{2}f=0$ for every cochain $f$. A $k$-cochain $f$ is a 
\textit{coboundary} if there exists a cochain $g$ such that $f=\partial g$;
a cochain $f$ is a \textit{cocycle }if\textit{\ }$\partial f=0$. We thus
have that every coboundary is a cocycle.

\section{The Keller-Maslov index\label{km}}

In this section we recall the definition of the Keller-Maslov index (seed 
\cite{Arnold,Leray,Souriau2} for a thorough study).

Using the ideas of Keller \cite{Keller}, Maslov introduced in \cite{Maslov}
an index, whose definition was clarified by \cite{Arnold} and extended by
Leray \cite{Leray5,Leray} and one of the authors \cite{JMPA}. (Dazord \cite%
{Dazord} has given similar results in the more general case of symplectic
bundles.) The Keller-Maslov index is designed to count the number of
intersections of a Lagrangian loop $\gamma $ with the Maslov cycle $\Sigma $.%
\textit{\ }It is defined in the following way. Let $W(n,\mathbb{C})$ be the
submanifold of $U(n,\mathbb{C})$ consisting of symmetric matrices: 
\begin{equation*}
W(n,\mathbb{C})=\left\{ u\in U(n,\mathbb{C}):u=u^{t}\right\} \text{.}
\end{equation*}%
($u^{t}=(u^{\ast })^{-1}$ the transpose of $u$). The mapping 
\begin{equation*}
\Lambda (n)\ni \ell =uX^{\ast }\longmapsto uu^{t}\in U(n,\mathbb{C})
\end{equation*}%
is a homeomorphism 
\begin{equation*}
\Lambda (n)\ni \ell =uX^{\ast }\longmapsto uu^{t}\in W(n,\mathbb{C})
\end{equation*}%
The composition of the isomorphism $\pi _{1}(\Lambda (n))\cong \pi _{1}(W(n,%
\mathbb{C}))$ induced by this homeomorphism with the isomorphism 
\begin{equation}
\pi _{1}(W(n,\mathbb{C}))\ni \lbrack \gamma ]\longmapsto \frac{1}{2\pi i}%
\oint_{\gamma }\frac{d(\det w)}{\det w}\in \mathbb{Z}  \label{isom1}
\end{equation}%
is an isomorphism 
\begin{equation*}
\func{ind}:\pi _{1}(\Lambda (n))\ni \lbrack \gamma ]\overset{\cong }{%
\longmapsto }\func{ind}[\gamma ]\in (\mathbb{Z},\mathbb{+)}\text{.}
\end{equation*}%
By definition, The Keller-Maslov index on $\Lambda (n)$ is the mapping $m$
that to each Lagrangian loop $\gamma $ associate the integer: 
\begin{equation}
m(\gamma )=\func{ind}[\gamma ]=\frac{1}{2\pi i}\oint_{\gamma }\frac{d(\det w)%
}{\det w}\text{.}  \label{indice}
\end{equation}

We can easily show that $m(\gamma )$ is zero in the strata of $\Lambda (n)$: 
\begin{equation*}
\func{Im}\gamma \subset \Lambda _{\ell }(n;k)\Longrightarrow m(\gamma )=0%
\text{.}
\end{equation*}

Let $\beta $ be the generator of $\pi _{1}(\Lambda (n))$ whose natural image
in $\mathbb{Z}$ is $+1$. Its Maslov index is 
\begin{equation}
m(\beta )=1\text{.}  \label{mabetaun}
\end{equation}
The Keller-Maslov index is an homotopic invariant and in particular $%
m(\gamma )=0$ if $\gamma $ is contractible to a point. Furthermore, this
index possesses the following additivity property: 
\begin{equation}
m(\gamma \ast \gamma ^{\prime })=m(\gamma )+m(\gamma ^{\prime })
\label{maprop}
\end{equation}
for all consecutive loops $\gamma $ and $\gamma ^{\prime }$ ($\gamma \ast
\gamma ^{\prime }$ being the loop $\gamma $ followed by the loop $\gamma
^{\prime }$).

One should note the difference between the Keller-Maslov index $m$ on $%
\Lambda (n)$ and the index $m_{V}$ of loop on a Lagrangian submanifold $V$
in $X\times X^{\ast }$. The index $m_{V}$ is defined in the following way:
Let $\ell (\cdot )$ be the continuous mapping that to each $z\in V$
associates the tangent space at $z$ : $\ell (z)=T_{z}V$. This mapping
associated to each loop $\gamma _{V}$ in $V$ a loop $\gamma $ in $\Lambda
(n) $, and by definition: 
\begin{equation}
m_{V}(\gamma _{V})=m(\gamma )\text{.}  \label{variete}
\end{equation}%
Souriau proved in \cite{Souriau3} that 
\begin{equation}
m_{V}(\gamma _{V})\text{ \textit{is an even integer when the manifold} }V%
\text{ \textit{is oriented}.}  \label{sousou}
\end{equation}%
In fact, we have the more precise result (see de Gosson \cite{Wiley}). In
what follows we denote by $\Lambda _{2q}(n)$ the connected $2q$-fold
covering of $\Lambda (n)$, $q=1,2,...$).

\begin{proposition}
Suppose that the mapping $\ell (\cdot ):V\longrightarrow \Lambda (n)$ is
lifted to a continuous mapping $\ell _{2q}(\cdot ):V\longrightarrow \Lambda
_{2q}(n).$ Then 
\begin{equation*}
m_{V}(\gamma _{V})\equiv 0\text{ \ }\func{mod}2q
\end{equation*}
for every loop $\gamma _{V}$ in $V$.
\end{proposition}

Note that if $\ell (\cdot )$ is lifted to a continuous mapping $\ell
_{\infty }(\cdot ):V\longrightarrow \Lambda _{\infty }(n)$ we have $m(\gamma
_{V})=0$ for every loop $\gamma _{V}$ in $V$ since the universal covering $%
\Lambda _{\infty }(n)$ is simply connected.

\section{Definition of the indices $\protect\mu _{\Lambda }$ and $\protect%
\mu _{Sp}$}

Let $\mathcal{C}(\Lambda (n))$ be the set of continuous paths $%
[0,1]\longrightarrow \Lambda (n)$. If $\lambda $ and $\lambda ^{\prime }$
are two consecutive paths (i.e., if $\lambda (1)=\lambda ^{\prime }(0)$) we
shall denote $\lambda \ast \lambda ^{\prime }$ the path $\lambda $ followed
by the path $\lambda ^{\prime }$: 
\begin{equation*}
\lambda \ast \lambda ^{\prime }(t)=\left\{ 
\begin{array}{c}
\lambda (2t)\text{ \ \textit{if} \ }0\leq t\leq \tfrac{1}{2}\medskip \\ 
\lambda ^{\prime }(2t-1)\text{ \ \textit{if} \ }\tfrac{1}{2}\leq t\leq 1%
\text{.}%
\end{array}%
\right.
\end{equation*}%
We denote by $\lambda ^{-1}$ the inverse of the path $\lambda $ : $\lambda
^{-1}(t)=\lambda (1-t)$. Finally, we shall write $\lambda \sim \lambda
^{\prime }$ when the paths $\lambda $ and $\lambda ^{\prime }$ are homotopic
with fixed endpoints.

\subsection{Axioms for the indices $\protect\mu _{\Lambda }$}

\begin{definition}
\label{definition1}A ''Lagrangian intersection index''\ is a mapping 
\begin{equation*}
\mu _{\Lambda }:\mathcal{C}(\Lambda (n))\times \Lambda (n)\ni (\lambda ,\ell
)\longmapsto \mu _{\Lambda (n)}(\lambda ,\ell )\in \mathbb{Z}
\end{equation*}
having the following four properties:

\begin{description}
\item[(L$_{1}$)] homotopic invariance\textbf{:}\ \textit{If the paths }$%
\lambda $\textit{\ and }$\lambda ^{\prime }$\textit{\ in }$\Lambda (n)$%
\textit{\ have the same endpoints, then} $\mu _{\Lambda }(\lambda ,\ell
)=\mu _{\Lambda }(\lambda ^{\prime },\ell )$ if and only if $\lambda \sim
\lambda ^{\prime }$;

\item[(L$_{2}$)] Additivity under composition: If $\lambda $ and $\lambda
^{\prime }$ are two consecutive paths, then for all $\ell \in \Lambda (n)$: 
\begin{equation*}
\mu _{\Lambda }(\lambda \ast \lambda ^{\prime },\ell )=\mu _{\Lambda
}(\lambda ,\ell )+\mu _{\Lambda }(\lambda ^{\prime },\ell )
\end{equation*}

\item[(L$_{3}$)] Zero in the stratum. If the path $\lambda $ remains in the
same stratum $\Lambda _{\ell }(n;k)$, then $\mu _{\Lambda }(\lambda ,\ell )$
is zero, i.e., 
\begin{equation*}
\dim \lambda (t)\cap \ell =k\text{ \ }(0\leq t\leq 1)\Longrightarrow \mu
_{\Lambda }(\lambda ,\ell )=0\text{.}
\end{equation*}

\item[(L$_{4}$)] Restriction to loops. If $\gamma $ is a loop, then $\mu
_{\Lambda }(\gamma ,\ell )=2m(\gamma )$ ($m(\gamma )$ the Maslov index of $%
\gamma $) for all $\ell $.\medskip

(We shall see in \ref{RSM} that the condition $\mu _{\Lambda }(\gamma ,\ell
)=2m(\gamma )$, and not $\mu _{\Lambda }(\gamma ,\ell )=m(\gamma )$, is
necessary for an index satisfying the axioms (L$_{1}$--L$_{3}$) to be an
integer)
\end{description}
\end{definition}

We note in particular that the axioms (L$_{2}$) and (L$_{4}$) imply that the
index is antisymmetric, i.e. , 
\begin{equation}
\mu _{\Lambda }(\lambda ^{-1},\ell )=-\mu _{\Lambda }(\lambda ,\ell )\text{.}
\label{mopp}
\end{equation}
Indeed, by (L$_{2}$) we have 
\begin{equation*}
\mu _{\Lambda (n)}(\lambda \ast \lambda ^{-1},\ell )=\mu _{\Lambda
(n)}(\lambda ,\ell )+\mu _{\Lambda (n)}(\lambda ^{-1},\ell )
\end{equation*}
but since the loop $\lambda \ast \lambda ^{-1}=\gamma $ is homotopic to a
point the axiom (L$_{4}$) implies that 
\begin{equation*}
\mu _{\Lambda (n)}(\lambda \ast \lambda ^{-1},\ell )=2m(\gamma )=0\text{.}
\end{equation*}

The system of axioms (L$_{1}$)--(L$_{4}$) is in fact equivalent to the
system of axioms obtained by replacing (L$_{1}$) by an apparently stronger
condition (L$_{1}^{\prime }$) we are going to describe below. Let us first
define the notion of ``homotopy in strata'':

\begin{definition}
Two Lagrangian paths $\lambda $ and $\lambda ^{\prime }$ are said to be
``homotopic in the strata relative to $\ell $'' (denoted $\lambda \approx
_{\ell }\lambda ^{\prime }$) if there exits a continuous mapping $%
h:[0,1]\times \lbrack 0,1]\longrightarrow \Lambda (n)$ such that 
\begin{equation*}
h(t,0)=\lambda (t)\text{ \ , }h(t,1)=\lambda ^{\prime }(t)\text{ \ \textit{%
for} \ }0\leq t\leq 1
\end{equation*}
and two integers $k_{0},k_{1}$ ($0\leq k_{0},k_{1}\leq n$) such that 
\begin{equation*}
h(0,s)\in \Lambda _{\ell }(n;k_{0})\text{ \ \textit{and} \ }h(1,s)\in
\Lambda _{\ell }(n;k_{1})\text{ \ \textit{for} \ }0\leq s\leq 1\text{.}
\end{equation*}
\end{definition}

The intersection indices $\mu _{\Lambda }$ have the following property that
makes (L$_{1}$) more precise:\medskip

\begin{proposition}
\begin{description}
\item 

\item[(L$_{1}^{\prime }$)] \textit{If the paths }$\lambda $\textit{\ and }$%
\lambda ^{\prime }$\textit{\ are homotopic in strata relative to }$\ell $%
\textit{, then} $\mu _{\Lambda }(\lambda ,\ell )=\mu _{\Lambda }(\lambda
^{\prime },\ell )$, i.e. , 
\begin{equation}
\lambda \approx _{\ell }\lambda ^{\prime }\Longrightarrow \mu _{\Lambda
}(\lambda ,\ell )=\mu _{\Lambda }(\lambda ^{\prime },\ell )\text{.}
\label{precis}
\end{equation}
\end{description}
\end{proposition}

\begin{proof}
Suppose that $\lambda \approx _{\ell }\lambda ^{\prime }$ and define the
paths $\varepsilon _{0}$ and $\varepsilon _{1}$ joining $\lambda ^{\prime
}(0)$ to $\lambda (0)$ and $\lambda (1)$ to $\lambda ^{\prime }(1)$,
respectively, by $\varepsilon _{0}(s)=h(0,1-s)$ and $\varepsilon
_{1}(s)=h(1,s)$ ($0\leq s\leq 1$). Then $\lambda \ast \varepsilon _{1}\ast
\lambda ^{\prime -1}\ast \varepsilon _{0}$ is homotopic to a point, and
hence, in view of (L$_{2}$) and (L$_{4}$): 
\begin{equation*}
\mu _{\Lambda }(\lambda ,\ell )+\mu _{\Lambda }(\varepsilon _{1},\ell )+\mu
_{\Lambda }(\lambda ^{\prime -1},\ell )+\mu _{\Lambda }(\varepsilon
_{0},\ell )=0\text{.}
\end{equation*}
But, in view of (L$_{3}$) 
\begin{equation*}
\mu _{\Lambda }(\varepsilon _{1},\ell )=\mu _{\Lambda }(\varepsilon
_{0},\ell )=0
\end{equation*}
and thus 
\begin{equation*}
\mu _{\Lambda }(\lambda ,\ell )+\mu _{\Lambda }(\lambda ^{\prime -1},\ell )=0
\end{equation*}
hence $\mu _{\Lambda }(\lambda ,\ell )=\mu _{\Lambda }(\lambda ^{\prime
},\ell )$\ in view of (\ref{mopp}).
\end{proof}

\subsection{A property of relative uniqueness}

The axioms (L$_{1}$--L$_{4}$) do not guaranty the uniqueness of an
intersection index. Indeed, if $\mu _{\Lambda }$ is such an index, then the
function $\mu _{\Lambda (n)}^{\prime }:\mathcal{C}(\Lambda (n))\times
\Lambda (n)\longrightarrow \mathbb{Z}$ defined by 
\begin{equation*}
\mu _{\Lambda }^{\prime }(\lambda ,\ell )=\mu (\lambda ,\ell )+\dim (\lambda
(1)\cap \ell )-\dim (\lambda (0)\cap \ell )
\end{equation*}%
is also an intersection index. We will prove in this section a result of
\textquotedblleft relative uniqueness\textquotedblright\ modulo a $1$%
-cochain.

Let us introduce the following notations: for every pair of points $(\ell
_{i},\ell _{j})$ of $\Lambda (n)$ we denote $\lambda _{ij}$ or $\lambda
_{ij}^{\prime }$ an arbitrary element of $\mathcal{C}(\Lambda (n))$ joining $%
\ell _{i}$ to $\ell _{j}$. The opposite path $\lambda _{ij}^{-1}$ of $%
\lambda _{ij}$ will then be denoted $\lambda _{ji}$. Using these notations
we can formulate a theorem relating all intersection indices to one another.
That is, show that there exists some function allowing us to relate one
intersection index to another, thereby giving us a relative uniqueness
theorem for intersection indices.

\begin{theorem}
\label{un}Let $\mu _{\Lambda }$ and $\mu _{\Lambda }^{\prime }$ be two
intersection indices on $\Lambda (n)$. For every $\ell \in \Lambda (n)$
there exists a function $f:\mathbb{N}\longrightarrow \mathbb{Z}$ such that: 
\begin{equation}
\mu _{\Lambda }(\lambda _{01},\ell )-\mu _{\Lambda }^{\prime }(\lambda
_{01},\ell )=f(\dim (\ell _{0}\cap \ell ))-f(\dim (\ell _{1}\cap \ell ))%
\text{.}  \label{kiki}
\end{equation}
\end{theorem}

To prove this fact we have to make sure that the left-hand side of (\ref%
{kiki}) does not depend on anything else than $\ell _{0},\ell _{1}$ (and of
course $\ell $). We have the essential

\begin{lemma}
for $\ell $ fixed, the difference $\mu _{\Lambda }(\lambda _{01},\ell )-\mu
_{\Lambda }^{\prime }(\lambda _{01},\ell )$ depends only on $(\ell _{0},\ell
_{1})$, and the $1$-cochain $\delta _{\ell }$ on $\Lambda (n)$ given by: 
\begin{equation}
\delta _{\ell }(\ell _{0},\ell _{1})=\mu _{\Lambda }(\lambda _{01},\ell
)-\mu _{\Lambda }^{\prime }(\lambda _{01},\ell )  \label{delta}
\end{equation}%
is an antisymmetric cocycle: 
\begin{equation}
\partial \delta _{\ell }(\ell _{0},\ell _{1},\ell _{2})=0\text{ \ and \ }%
\delta _{\ell }(\ell _{0},\ell _{1})=-\delta _{\ell }(\ell _{1},\ell _{0})
\label{cocasse}
\end{equation}%
for all $(\ell _{0},\ell _{1},\ell _{2})\in \Lambda (n)^{3}$.
\end{lemma}

\begin{proof}[Proof of the lemma]
\smallskip Suppose that we replace the path $\lambda _{01}$ with another
path $\lambda _{01}^{\prime }$ joining $\ell _{0}$ to $\ell _{1}$. Then,
there exists a loop $\gamma _{1}$ passing by $\ell _{1}$ and such that $%
\lambda _{01}^{\prime }\sim \lambda _{01}\ast \gamma _{1}$ and we have, by
the axioms (L$_{1}$) and (L$_{4}$): 
\begin{equation*}
\left\{ 
\begin{array}{c}
\mu _{\Lambda }(\lambda _{01}^{\prime },\ell )=\mu _{\Lambda }(\lambda
_{01},\ell )+2m(\gamma _{1})\medskip \\ 
\mu _{\Lambda }^{\prime }(\lambda _{01}^{\prime },\ell )=\mu _{\Lambda
}^{\prime }(\lambda _{01},\ell )+2m(\gamma _{1})%
\end{array}%
\right. .
\end{equation*}%
Hence, the difference 
\begin{equation*}
\mu _{\Lambda }(\lambda _{01}^{\prime },\ell )-\mu _{\Lambda }^{\prime
}(\lambda _{01}^{\prime },\ell )=\mu _{\Lambda }(\lambda _{01},\ell )-\mu
_{\Lambda }^{\prime }(\lambda _{01},\ell )=\delta _{\ell }(\ell _{0},\ell
_{1})
\end{equation*}%
depends only on $\ell _{0}$ and $\ell _{1}$ (and $\ell $). We now show that
the coboundary of that $1$-cochain is zero: $\partial \delta _{\ell }=0$.
Let $(\ell _{0},\ell _{1},\ell _{2})$ be an arbitrary triple of points in $%
\Lambda (n)$. By the properties (L$_{2}$) and (L$_{4}$) we see that 
\begin{equation}
\left\{ 
\begin{array}{c}
\mu _{\Lambda }(\lambda _{01},\ell )-\mu _{\Lambda }(\lambda _{02},\ell
)+\mu _{\Lambda }(\lambda _{12},\ell )=2m(\gamma _{0})\medskip \\ 
\mu _{\Lambda }^{\prime }(\lambda _{01},\ell )-\mu _{\Lambda }^{\prime
}(\lambda _{02},\ell )+\mu _{\Lambda }^{\prime }(\lambda _{12},\ell
)=2m(\gamma _{0})%
\end{array}%
\right.  \label{accol}
\end{equation}%
where $\gamma _{0}=\lambda _{01}\ast \lambda _{12}\ast \lambda _{20}$ is a
loop with origin $\ell _{0}$. Substracting the second equality (\ref{accol})
from the first and using the definition of $\delta _{\ell }$ we obtain: 
\begin{equation*}
\partial \delta _{\ell }(\ell _{0},\ell _{1},\ell _{2})=0\text{.}
\end{equation*}%
The antisymmetry of $\delta _{\ell }$ is an immediate consequence of (\ref%
{mopp}).
\end{proof}

\begin{proof}[Proof of the theorem]
(1) By the definition of $\delta _{\ell }$ we have 
\begin{equation*}
\mu _{\Lambda }(\lambda _{01},\ell )-\mu _{\Lambda }(\lambda _{01}^{\prime
},\ell )=\delta _{\ell }(\ell _{0},\ell _{1})\text{.}
\end{equation*}
Since $\delta _{\ell }$ is a cocycle this equality can be written as 
\begin{equation*}
\mu _{\Lambda }(\lambda _{01},\ell )-\mu _{\Lambda }(\lambda _{01}^{\prime
},\ell )=\delta _{\ell }(\ell _{0},\ell )-\delta _{\ell }(\ell _{1},\ell )%
\text{.}
\end{equation*}
We shall show that the $0$-cochain $\ell _{0}\longmapsto \delta _{\ell
}(\ell _{0},\ell )$ is constant in each stratum $\Lambda _{\ell }(n;k)$ and
the theorem will follow. We see, by the definition of $\delta _{\ell }$ that 
\begin{equation*}
\delta _{\ell }(\ell _{0},\ell )=\mu _{\Lambda }(\lambda ,\ell )-\mu
_{\Lambda }^{\prime }(\lambda ,\ell )
\end{equation*}
where $\lambda $ is \ an arbitrary path joining $\ell _{0}$ to $\ell $ in $%
\Lambda (n)$. Let $\lambda ^{\prime }$ be a path joining $\ell _{0}$ to $%
\ell _{0}^{\prime }$ in $\Lambda _{\ell }(n;k)$. Then 
\begin{equation*}
\delta _{\ell }(\ell _{0}^{\prime },\ell )=\mu _{\Lambda }(\lambda ^{\prime
}\ast \lambda ,\ell )-\mu _{\Lambda }^{\prime }(\lambda ^{\prime }\ast
\lambda ,\ell )
\end{equation*}
and in view of the additivity property (L$_{2}$) 
\begin{equation*}
\delta _{\ell }(\ell _{0}^{\prime },\ell )=\mu _{\Lambda }(\lambda ^{\prime
},\ell )-\mu _{\Lambda }^{\prime }(\lambda ^{\prime },\ell )+\delta _{\ell
}(\ell _{0},\ell )\text{.}
\end{equation*}
But, by the axiom (L$_{3}$) 
\begin{equation*}
\mu _{\Lambda }(\lambda ^{\prime },\ell )=\mu _{\Lambda }^{\prime }(\lambda
^{\prime },\ell )=0
\end{equation*}
and we thus have $\delta _{\ell }(\ell _{0}^{\prime },\ell )=\delta _{\ell
}(\ell _{0},\ell )$, from which we conclude that the mapping $\ell
_{0}\longmapsto \delta _{\ell }(\ell _{0},\ell )$ is constant in the stratum 
$\Lambda _{\ell }(n;k)$.
\end{proof}

\subsection{The signature cocycle\label{cosign}}

We now present an interesting result that will generalize the axiom (L$_{1}$%
) and that will enable us to give a notion of signature for every triple of
elements of elements in $\Lambda (n)$.

\begin{lemma}
\label{3chemins}Let $\ell _{0}$, $\ell _{1}$, $\ell _{2}$ be three elements
of $\Lambda (n)$. Suppose that the paths $\lambda _{ij}$ and $\lambda
_{ij}^{\prime }$ ($0\leq i,j\leq 2$) are such that 
\begin{equation}
\lambda _{01}\ast \lambda _{12}\ast \lambda _{20}\sim \lambda _{01}^{\prime
}\ast \lambda _{12}^{\prime }\ast \lambda _{20}^{\prime }\text{.}
\label{etoile}
\end{equation}
Then both sums 
\begin{equation*}
\Sigma =\mu _{\Lambda }(\lambda _{01},\ell _{2})+\mu _{\Lambda }(\lambda
_{12},\ell _{0})+\mu _{\Lambda }(\lambda _{20},\ell _{1})
\end{equation*}
and 
\begin{equation*}
\Sigma ^{\prime }=\mu _{\Lambda }(\lambda _{01}^{\prime },\ell _{2})+\mu
_{\Lambda }(\lambda _{12}^{\prime },\ell _{0})+\mu _{\Lambda }(\lambda
_{20}^{\prime },\ell _{1})
\end{equation*}
are equal: $\Sigma =\Sigma ^{\prime }$. Furthermore, $\Sigma $ does not
depend on the choice of $\mu _{\Lambda }$.
\end{lemma}

\begin{proof}
There exist loops $\gamma _{0}^{1}$, $\gamma _{0}^{2}$ and $\gamma _{0}^{3}$
passing through $\ell _{0}$ such that 
\begin{equation*}
\left\{ 
\begin{array}{c}
\lambda _{01}^{\prime }\sim \gamma _{0}^{1}\ast \lambda _{01}\medskip \\ 
\lambda _{12}^{\prime }\sim \lambda _{10}\ast \gamma _{0}^{2}\ast \lambda
_{01}\ast \lambda _{12}\medskip \\ 
\lambda _{20}^{\prime }\sim \lambda _{20}\ast \gamma _{0}^{3}%
\end{array}
\right. \text{ .}
\end{equation*}
In view of (\ref{etoile}) and the properties (L$_{1}$) and (L$_{4}$) we
obtain 
\begin{align*}
\Sigma ^{\prime }=\Sigma +2(m(\gamma _{0}^{1})+m(\gamma _{0}^{2})+m(\gamma
_{0}^{3})) \\
=\Sigma +2m(\gamma _{0}^{1}\ast \gamma _{0}^{2}\ast \gamma _{0}^{3})\text{.}
\end{align*}
But (\ref{etoile}) implies that $\gamma _{0}^{1}\ast \gamma _{0}^{2}\ast
\gamma _{0}^{3}\sim 0$ and hence that $\mu _{\Lambda (n)}(\gamma
_{0}^{1}\ast \gamma _{0}^{2}\ast \gamma _{0}^{3})=0$ and $\Sigma ^{\prime
}=\Sigma $. It remains to prove that $\Sigma $ is independent of the choice
of index $\mu _{\Lambda (n)}$. If $\mu _{\Lambda (n)}^{\prime }$ is another
intersection index\ and $\Sigma ^{\prime }$\ the associated sum the relative
uniqueness theorem \ref{un} tells us that 
\begin{multline*}
\Sigma -\Sigma ^{\prime }=f(\dim (\ell _{0}\cap \ell _{2}))-f(\dim (\ell
_{1}\cap \ell _{2}))+f(\dim (\ell _{1}\cap \ell _{0})) \\
-f(\dim (\ell _{2}\cap \ell _{0}))+f(\dim (\ell _{2}\cap \ell _{1}))-f(\dim
(\ell _{0}\cap \ell _{1}))
\end{multline*}
that is: $\Sigma -\Sigma ^{\prime }=0$.
\end{proof}

This result motivates the following definition

\begin{definition}
Let $(\ell _{0},\ell _{1},\ell _{2})\in \Lambda (n)^{3}$ and $\lambda _{01}$%
, $\lambda _{12}$, $\lambda _{20}$ be elements of $\mathcal{C}$ $(\Lambda
(n))$ such that $\lambda _{01}\ast \lambda _{12}\ast \lambda _{20}$ be
homotopic to a point. The sum 
\begin{equation}
\func{sign}(\ell _{0},\ell _{1},\ell _{2})=\mu _{\Lambda }(\lambda
_{01},\ell _{2})+\mu _{\Lambda }(\lambda _{12},\ell _{0})+\mu _{\Lambda
}(\lambda _{20},\ell _{1})  \label{siglag}
\end{equation}
is called the signature of the triple $(\ell _{0},\ell _{1},\ell _{2})$.
This signature is independent of the choice of intersection index $\mu
_{\Lambda }$ on $\Lambda (n)$.
\end{definition}

\begin{remark}
In Section \ref{Kashiwa} we will show that $\func{sign}$ is just the
Demazure--Kashiwara signature $\tau $ (\cite{Marle,LV}) arising in
symplectic geometry.
\end{remark}

By (\ref{mopp}), the signature is an antisymmetric $2$-cochain 
\begin{equation*}
\func{sign}(\ell _{\varepsilon (0)},\ell _{\varepsilon (1)},\ell
_{\varepsilon (2)})=(-1)^{sgn(\varepsilon )}\func{sign}(\ell _{0},\ell
_{1},\ell _{2})
\end{equation*}
for every permutation $\varepsilon $ of $\left\{ 0,1,2\right\} $.
Furthermore, it possesses the following essential property

\begin{proposition}
\label{propos}The signature is a $2$-cocycle: 
\begin{equation}
\partial \func{sign}(\ell _{0},\ell _{1},\ell _{2},\ell _{3})=0
\label{sicoc}
\end{equation}
that is 
\begin{multline*}
\func{sign}(\ell _{1},\ell _{2},\ell _{3})-\func{sign}(\ell _{0},\ell
_{2},\ell _{3}) \\
+\func{sign}(\ell _{0},\ell _{1},\ell _{3})-\func{sign}(\ell _{0},\ell
_{1},\ell _{2})=0\text{.}
\end{multline*}
\end{proposition}

\begin{proof}
By the definition of the coboundary operator $\partial $ we have 
\begin{multline*}
\partial \func{sign}(\ell _{0},\ell _{1},\ell _{2},\ell _{3})=\func{sign}%
(\ell _{1},\ell _{2},\ell _{3})-\func{sign}(\ell _{0},\ell _{2},\ell _{3}) \\
+\func{sign}(\ell _{0},\ell _{1},\ell _{3})-\func{sign}(\ell _{0},\ell
_{1},\ell _{2})
\end{multline*}
and a calculation making use of (\ref{mopp}) and the axiom (L$_{2}$) gives
us 
\begin{multline*}
\partial \func{sign}(\ell _{0},\ell _{1},\ell _{2},\ell _{3})= \\
\mu _{\Lambda (n)}(\gamma _{1},\ell _{3})+\mu _{\Lambda (n)}(\gamma
_{2},\ell _{1})+\mu _{\Lambda (n)}(\gamma _{3},\ell _{2})+\mu _{\Lambda
(n)}(\gamma _{4},\ell _{0})
\end{multline*}
where $\gamma _{1}=\gamma _{12}\ast \gamma _{20}\ast \gamma _{01}$, $\gamma
_{2}=\gamma _{23}\ast \gamma _{30}\ast \gamma _{02}$, $\gamma _{3}=\gamma
_{31}\ast \gamma _{10}\ast \gamma _{03}$ and $\gamma _{4}=\gamma _{32}\ast
\gamma _{21}\ast \gamma _{13}$ are loop, hence by (L$_{4}$) 
\begin{equation*}
\partial \func{sign}(\ell _{0},\ell _{1},\ell _{2},\ell _{3})=2(m(\gamma
_{1})+m(\gamma _{2})+m(\gamma _{3})+m(\gamma _{3}))\text{.}
\end{equation*}
The loops $\gamma _{1}$, $\gamma _{2}$, $\gamma _{3}$ and $\gamma _{4}$
being contractible to a point, the Maslov indices of this loops are zero and
thereby $\partial \func{sign}=0$.
\end{proof}

\subsection{The indices $\protect\mu _{Sp}$}

We shall denote by $\mathcal{C}(Sp(n))$ the set of continuous paths $%
[0,1]\longrightarrow Sp(n)$.

\begin{definition}
\label{definition2} A symplectic intersection index is a mapping 
\begin{equation*}
\mu _{Sp}:\mathcal{C}(Sp(n))\times \Lambda (n)\ni (\sigma ,\ell )\longmapsto
\mu _{Sp}(\lambda ,\ell )\in \mathbb{Z}
\end{equation*}%
satisfying the following four axioms:

\begin{description}
\item[(S$_{1}$)] Homotopic invariance\textbf{.} If the symplectic paths $%
\sigma $ and $\sigma ^{\prime }$ are homotopic with fixed endpoints, then $%
\mu _{Sp}(\sigma ,\ell )=\mu _{Sp}(\sigma ^{\prime },\ell )$ for all $\ell
\in \Lambda (n)$.

\item[(S$_{2}$)] Additivity by concatenation\textbf{.} If $\sigma $ and $%
\sigma ^{\prime }$ are two consecutive symplectic paths, and if $\sigma \ast
\sigma ^{\prime }$ is the path $\sigma $ followed by the path $\sigma
^{\prime }$ then 
\begin{equation*}
\mu _{Sp}(\sigma \ast \sigma ^{\prime },\ell )=\mu _{Sp}(\sigma ,\ell )+\mu
_{Sp}(\sigma ^{\prime },\ell )
\end{equation*}
for all $\ell \in \Lambda (n)$.

\item[(S$_{3}$)] Zero in the stratum\textbf{. }If $\sigma $ and $\ell $ are
such that $\func{Im}(\sigma \ell )\subset \Lambda (n)_{\ell }$, then $\mu
_{Sp}(\sigma ,\ell )=0$.

\item[(S$_{4}$)] Restriction to loops. If $\psi $ is a loop in $Sp$, then $%
\mu _{Sp}(\psi ,\ell )=2m(\psi \ell )$ for all $\ell \in \Lambda (n)$%
.\medskip
\end{description}
\end{definition}

As the homotopy axiom (L$_{1}$) for the Lagrangian indices that may be
replaced by an axiom (L$_{1}^{\prime }$) of homotopy in the strata the axiom
(S$_{1}$) above may be strengthened by the following statement

\begin{description}
\item[(S$_{1}^{\prime }$)] \textit{If the symplectic paths }$\sigma $\textit{%
\ and }$\sigma ^{\prime }$\textit{\ are such that }$\sigma \ell $ and $%
\sigma ^{\prime }\ell $ \textit{are homotopic in strata relative to }$\ell $%
\textit{, then }$\mu _{Sp}(\sigma ,\ell )=\mu _{Sp}(\sigma ^{\prime },\ell )$%
.
\end{description}

The proof of this is identical to that of (L$_{1}^{\prime }$). What is
noteworthy is that the data of an intersection index on $\Lambda (n)$ is
equivalent to that of an intersection index on $Sp(n)$. Indeed, let $\mu
_{\Lambda }$ be an intersection index on $\Lambda (n)$ and let $\sigma \in 
\mathcal{C}(Sp(n))$ be a symplectic path. Then then function 
\begin{equation}
\mathcal{C}(Sp(n))\times \Lambda (n)\ni (\sigma ,\ell )\longmapsto \mu
(\sigma \ell ,\ell )\in \mathbb{Z}  \label{formuleun}
\end{equation}%
($\sigma \ell $ being the path $t\longmapsto \sigma (t)\ell $) is an
intersection index on $Sp(n)$.

Conversely, to each intersection index $\mu _{Sp}$ we may associate an
intersection index $\mu _{\Lambda (n)}$ on $\Lambda (n)$ in the following
way. To each $\ell \in \Lambda (n)$ we have a fibration 
\begin{equation*}
Sp(n)\longrightarrow Sp/St(\ell )=\Lambda (n)
\end{equation*}
($St(\ell )$ being the stabilizer of $\ell $ in $Sp$), hence, every path $%
\lambda \in \mathcal{C}(\Lambda (n))$ can be lifted to a path $\sigma \in 
\mathcal{C}(Sp(n))$ such that $\lambda =\sigma \ell $. One verifies that the
mapping 
\begin{equation}
\mu _{\Lambda (n)}:\mathcal{C}(\Lambda (n))\times \Lambda (n)\ni (\lambda
,\ell )\longmapsto \mu _{Sp}(\sigma ,\ell )\in \mathbb{Z}
\label{formuledeux}
\end{equation}
defines an intersection index on the manifold $\Lambda (n)$.\medskip

As in the Lagrangian case, we have a relative uniqueness result for the
symplectic path intersection indices

\begin{corollary}
Let $\mu _{Sp}$ and $\mu _{Sp}^{\prime }$ be two symplectic intersection
indices. There exists a function $f$ such that 
\begin{multline}
\mu _{Sp}(\sigma ,\ell )=\mu _{Sp}^{\prime }(\sigma ,\ell )+  \label{21} \\
f(\dim (\sigma (0)\ell \cap \ell ))-f(\dim (\sigma (0)\ell \cap \ell ))\text{%
.}  \notag
\end{multline}
\end{corollary}

\begin{proof}
It is a direct consequence of the relative uniqueness theorem for Lagrangian
intersection indices.
\end{proof}

\section{The indices $\protect\mu _{\infty }$ on $\Lambda _{\infty }(n)$}

We shall in this section use the same notations as in the previous one for
paths joining some point $\ell _{i}$ to another point $\ell _{j}$ in $%
\Lambda (n)$. Recall that $\Lambda _{\infty }(n)$ denotes the universal
covering of $\Lambda (n)$. We are going to use the intersection index for
Lagrangian paths to define an index $\mu _{\infty }:\Lambda _{\infty
}(n)\times \Lambda _{\infty }(n)\longrightarrow \mathbb{Z}$.

\subsection{The definition of $\protect\mu _{\infty }$}

To each Lagrangian intersection index $\mu _{\Lambda }$ we may associate a $%
1 $-cochain $\mu _{\infty }$ on the universal covering $\pi :\Lambda
_{\infty }(n)\longrightarrow \Lambda (n)$ in the following way. We choose
once and for all a base point $\ell _{0}$ for $\Lambda (n)$. The elements of 
$\Lambda _{\infty }(n)$ are then homotopy classes with fixed endpoints of
continuous paths of origin $\ell _{0}$ in $\Lambda (n)$. The projection $\pi 
$ from $\Lambda _{\infty }(n)$ onto $\Lambda (n)$ is defined by $\pi (\ell
_{\infty })=\ell $ if $\ell _{\infty }$ is the homotopy class of a path $%
\lambda $ joining $\ell _{0}$ to $\ell $.

\begin{definition}
The Leray index associated to $\mu _{\Lambda }$ is the $1$-cochain $\mu
_{\infty }$ on $\Lambda _{\infty }(n)$ defined by 
\begin{equation}
\mu _{\infty }(\ell _{1,\infty },\ell _{2,\infty })=-\mu _{\Lambda
(n)}(\lambda _{12},\ell _{0})+\func{sign}(\ell _{0},\ell _{1},\ell _{2})
\label{oignon}
\end{equation}
where $\ell _{1,\infty }$ and $\ell _{2,\infty }$ are the homotopy classes
of the paths $\lambda _{01}$ and $\lambda _{02}$, respectively, and $\lambda
_{12}\sim \lambda _{10}\ast \lambda _{02}$. Also, $\func{sign}$ is the
signature cocycle associated to $\mu _{\Lambda }$, as defined in Subsection %
\ref{cosign}.
\end{definition}

Note that the Leray indices are antisymmetric 
\begin{equation}
\mu _{\infty }(\ell _{1,\infty },\ell _{2,\infty })=-\mu _{\infty }(\ell
_{2,\infty },\ell _{1,\infty })  \label{antsm}
\end{equation}
for all $(\ell _{1,\infty },\ell _{2,\infty })\in \Lambda (n)_{\infty }^{2}$
in view of (\ref{mopp}) and the fact that $\limfunc{sign}$ is antisymmetric.

By theorem \ref{un} and the fact that the definition of the signature
cocycle is intrinsic, two Leray indices $\mu _{\infty }$ and $\mu _{\infty
}^{\prime }$ associated to two intersection indices $\mu _{\Lambda }$ and $%
\mu _{\Lambda }^{\prime }$ are such that 
\begin{equation}
\mu _{\infty }(\ell _{1,\infty },\ell _{2,\infty })-\mu _{\infty }^{\prime
}(\ell _{1,\infty },\ell _{2,\infty })=f(\dim (\ell _{2}\cap \ell
_{0}))-f(\dim (\ell _{1}\cap \ell _{0}))  \notag
\end{equation}
for some function $f:\mathbb{N}\longrightarrow \mathbb{Z}$.\smallskip

The Leray indices have the following properties

\begin{proposition}
(1) The coboundary of $\mu _{\infty }$ descends to $\Lambda (n)$, i.e. , $%
\partial \mu _{\infty }=\pi ^{\ast }\func{sign}$. That is, 
\begin{equation}
\mu _{\infty }(\ell _{1,\infty },\ell _{2,\infty })-\mu _{\infty }(\ell
_{1,\infty },\ell _{3,\infty })+\mu _{\infty }(\ell _{2,\infty },\ell
_{3,\infty })=\func{sign}(\ell _{1},\ell _{2},\ell _{3})  \label{unideux}
\end{equation}
for all $(\ell _{1,\infty },\ell _{2,\infty },\ell _{3,\infty })\in \Lambda
_{\infty }(n)^{3}$;

(2) The action of $\mathbb{\pi }_{1}(\Lambda (n))$ on $\mu _{\infty }$ is
given by 
\begin{equation}
\mu _{\infty }(\gamma _{1}\ell _{1,\infty },\gamma _{2}\ell _{2,\infty
})=\mu _{\infty }(\ell _{1,\infty },\ell _{2,\infty })+2(m(\gamma
_{1})-m(\gamma _{2}))  \label{mabeta}
\end{equation}
for all $\gamma _{1},\gamma _{2}\in \mathbb{\pi }_{1}(\Lambda (n))$.
\end{proposition}

\begin{proof}
\textit{(1)} By the definition (\ref{oignon}) of $\mu _{\infty }$ we have 
\begin{multline*}
\mu _{\infty }(\ell _{1,\infty },\ell _{2,\infty })-\mu _{\infty }(\ell
_{1,\infty },\ell _{3,\infty })+\mu _{\infty }(\ell _{2,\infty },\ell
_{3,\infty })= \\
-(\mu _{\Lambda }(\lambda _{12},\ell _{0})-\mu _{\Lambda }(\lambda
_{13},\ell _{0})+\mu _{\Lambda }(\lambda _{23},\ell _{0}))+ \\
\func{sign}(\ell _{0},\ell _{2},\ell _{3})-\func{sign}(\ell _{0},\ell
_{1},\ell _{3})+\func{sign}(\ell _{0},\ell _{2},\ell _{3})\text{.}
\end{multline*}
Formula (\ref{unideux}) follows: by the axioms (L$_{2}$) and (L$_{4}$), and
taking into consideration the fact that $\gamma _{3}=\lambda _{31}\ast
\lambda _{12}\ast \lambda _{23}$ is homotopic to a point, we have 
\begin{equation*}
\mu _{\Lambda }(\lambda _{12},\ell _{0})-\mu _{\Lambda }(\lambda _{13},\ell
_{0})+\mu _{\Lambda }(\lambda _{23},\ell _{0})=2m(\gamma _{3})=0\text{.}
\end{equation*}
By the cocycle property of $\func{sign}$ this yields 
\begin{equation*}
\func{sign}(\ell _{0},\ell _{2},\ell _{3})-\func{sign}(\ell _{0},\ell
_{1},\ell _{3})+\func{sign}(\ell _{0},\ell _{2},\ell _{3})=\func{sign}(\ell
_{1},\ell _{2},\ell _{3})\text{.}
\end{equation*}
\textit{(2)} In view of the antisymmetry of $\mu _{\infty }$ it suffices to
show that 
\begin{equation}
\mu _{\infty }(\gamma \ell _{1,\infty },\ell _{2,\infty })=\mu _{\infty
}(\ell _{1,\infty },\ell _{2,\infty })+2m(\gamma )  \label{cp}
\end{equation}
for all $\gamma \in \pi _{1}(\Lambda (n))$. Let $\gamma $ be a loop with
origin $\ell _{0}$. If $\lambda _{01}$ is a representative of $\ell
_{1,\infty }$, then $\gamma \ast \lambda _{01}$is a representative of $%
\gamma \ell _{1,\infty }$ and, by the definition of $\mu _{\infty }$ : 
\begin{align*}
\mu _{\infty }(\gamma \ell _{1,\infty },\ell _{2,\infty })=-\mu _{\Lambda
}(\lambda _{10}\ast \gamma \ast \lambda _{02},\ell _{0})+\func{sign}(\ell
_{0},\ell _{1},\ell _{2}) \\
=-\mu _{\Lambda }(\lambda _{10}\ast \lambda _{02},\ell _{0})-2m(\gamma )+%
\func{sign}(\ell _{0},\ell _{1},\ell _{2}) \\
=\mu _{\infty }(\ell _{1,\infty },\ell _{2,\infty })+2m(\gamma )\text{.}
\end{align*}
Formula (\ref{cp}) then follows by the definition of the Maslov index.
\end{proof}

It is noteworthy that, conversely, the datum of a $1$-cochain $\mu _{\infty
} $ on $\Lambda _{\infty }(n)$ satisfying the property (\ref{mabeta}) above
together with a simple topological property, enables us to construct an
index, and by theorem \ref{un}, \emph{all} the Lagrangian intersection
indices (and hence all symplectic intersection indices).

\begin{theorem}
\label{deep}Let $\mu _{\infty }$ be a $1$-cochain on $\Lambda (n)$ having
the property (\ref{mabeta}), and being locally constant on each of the sets $%
\Lambda _{\infty }^{2}(n;k)$ ($0\leq k\leq n$) defined in Subsection \ref%
{not}. For $(\lambda _{12},\ell )\in \mathcal{C}(\Lambda (n))\times \Lambda
(n)$ let us define $\ell _{\infty }$, $\ell _{1,\infty }$ and $\ell
_{2,\infty }$ in the following way:

\begin{enumerate}
\item $\ell _{\infty }$ is an arbitrary element of $\Lambda _{\infty }(n)$
covering $\ell $.

\item $\ell _{1,\infty }$ is the equivalence class of an arbitrary path $%
\lambda _{01}\in $ $\mathcal{C}(\Lambda (n))$ joining $\ell _{0}$ to $\ell
_{1}$.

\item $\ell _{2,\infty }$ is the equivalence class of $\lambda _{02}=\lambda
_{01}\ast \lambda _{12\text{ }}$.
\end{enumerate}

Then the formula 
\begin{equation}
\mu _{\Lambda }(\lambda _{12},\ell )=\mu _{\infty }(\ell _{2,\infty },\ell
_{\infty })-\mu _{\infty }(\ell _{1,\infty },\ell _{\infty })
\label{foufoun}
\end{equation}
defines an intersection index. In particular we have $\partial \mu _{\infty
}=\pi ^{\ast }\func{sign}$.
\end{theorem}

\begin{proof}
That $\mu _{\Lambda (n)}(\lambda _{12},\ell )$ is independent of the choice
of the element $\ell _{\infty }$ of the universal covering $\Lambda _{\infty
}(n)$ covering $\ell $ follows immediately from (\ref{mabeta}). Indeed, if $%
\ell _{\infty }^{\prime }$ covers $\ell $ then there exists an element $%
\gamma \in \pi _{1}(\Lambda (n))$ such that $\ell _{\infty }^{\prime
}=\gamma \ell _{\infty }$ and hence 
\begin{equation*}
\left\{ 
\begin{array}{c}
\mu _{\infty }(\ell _{2,\infty },\ell _{\infty }^{\prime })=\mu _{\infty
}(\ell _{2,\infty },\ell _{\infty })-2m(\gamma )\smallskip \\ 
\mu _{\infty }(\ell _{1,\infty },\ell _{\infty }^{\prime })=\mu _{\infty
}(\ell _{1,\infty },\ell _{\infty })-2m(\gamma )%
\end{array}
\right. \text{\ }
\end{equation*}
and consequently 
\begin{equation*}
\mu _{\infty }(\ell _{2,\infty },\ell _{\infty }^{\prime })-\mu _{\infty
}(\ell _{1,\infty },\ell _{1,\infty }^{\prime })=\mu _{\infty }(\ell
_{2,\infty },\ell _{\infty })-\mu _{\infty }(\ell _{1,\infty },\ell _{\infty
})\text{.}
\end{equation*}
Let us show that $\mu _{\Lambda }(\lambda _{12},\ell )$ is also independent
of the choice of $\lambda _{01}$ and thereby of the choice of the element $%
\ell _{1,\infty }$ covering $\ell _{1}$. We replace $\lambda _{01}$ by some
path $\lambda _{01}^{\prime }$ and denote $\ell _{1,\infty }^{\prime }$ the
element of $\Lambda _{\infty }(n)$ it defines. Then, there exists an element 
$\gamma \ $\ of $\pi _{1}(\Lambda (n))$ such that $\ell _{1,\infty }^{\prime
}=\gamma \ell _{1,\infty }$; $\ell _{2,\infty }$ will thus be replaced by $%
\ell _{2,\infty }^{\prime }=\gamma \ell _{2,\infty }$. Using once again (\ref%
{mabeta}) we have 
\begin{equation*}
\mu _{\infty }(\ell _{2,\infty }^{\prime },\ell _{\infty })-\mu _{\infty
}(\ell _{1,\infty }^{\prime },\ell _{\infty })=\mu _{\infty }(\ell
_{2,\infty },\ell _{\infty })-\mu _{\infty }(\ell _{1,\infty },\ell _{\infty
})\text{.}
\end{equation*}
It remains to prove that the function $\mu _{\Lambda }$ defined by (\ref%
{foufoun}) satisfies the axioms (L$_{1}$)--(L$_{4}$).

\begin{enumerate}
\item[(L$_{1}$)] Let us replace the path $\lambda _{12}$ by any path $%
\lambda _{12}^{\prime }$ homotopic (with fixed endpoints) to $\lambda _{12}$%
. Then $\lambda _{02}=\lambda _{01}\ast \lambda _{12\text{ }}$is replaced by
a homotopic path $\lambda _{02}^{\prime }=\lambda _{01}\ast \lambda _{12%
\text{ }}^{\prime }$and the homotopy class $\ell _{2,\infty }$ does not
change. Consequently, $\mu _{\Lambda }(\lambda _{12}^{\prime },\ell )=\mu
_{\Lambda }(\lambda _{12},\ell )$.

\item[(L$_{2}$)] Consider two consecutive paths $\lambda _{12}$ and $\lambda
_{23}$. The index $\mu _{\Lambda }(\lambda _{23},\ell )$ is given by 
\begin{equation*}
\mu _{\Lambda }(\lambda _{23},\ell )=\mu _{\infty }(\ell _{3,\infty },\ell
_{\infty })-\mu _{\infty }(\ell _{2,\infty }^{\prime },\ell _{\infty })
\end{equation*}
where $\ell _{2,\infty }^{\prime }$ is the homotopy class of an arbitrary
path $\lambda _{02}^{\prime }$ and $\ell _{3,\infty }$ that of $\lambda
_{02}^{\prime }\ast \lambda _{23}$. Let us choose $\lambda _{02}^{\prime
}=\lambda _{02}$. Then $\ell _{2,\infty }^{\prime }=\ell _{2,\infty }$ and
we obtain 
\begin{multline*}
\mu _{\Lambda }(\lambda _{12},\ell )+\mu _{\Lambda }(\lambda _{23},\ell
)=\mu _{\infty }(\ell _{2,\infty },\ell _{\infty })-\mu _{\infty }(\ell
_{1,\infty },\ell _{\infty }) \\
+\mu _{\infty }(\ell _{3,\infty },\ell _{\infty })-\mu _{\infty }(\ell
_{2,\infty },\ell _{\infty })
\end{multline*}
that is, 
\begin{equation*}
\mu _{\Lambda }(\lambda _{12},\ell )+\mu _{\Lambda }(\lambda _{23},\ell
)=\mu _{\Lambda }(\lambda _{13},\ell )\text{.}
\end{equation*}

\item[(L$_{3}$)] Let $\lambda _{12}$ be a path in the strata $\Lambda _{\ell
}(n;k)$ and denote $\ell _{\infty }(t)$ the equivalence class of $\lambda
_{01}\ast \lambda _{12}(t)$. The mapping $t\longmapsto \ell _{\infty }(t)$
being continuous, the composition mapping $t\longmapsto \mu _{\infty }(\ell
_{\infty }(t),\ell _{\infty })$ is locally constant on the interval $[0,1]$.
It follows that it is constant on the same interval since $\Lambda _{\ell
}(n;k)$ is connected. Its values is 
\begin{equation*}
\mu _{\infty }(\ell _{\infty }(0),\ell _{\infty })=\mu _{\infty }(\ell
_{\infty }(1),\ell _{\infty })
\end{equation*}
whereof $\mu (\lambda _{12},\ell )=0$.

\item[(L$_{4}$)] \textit{\ }Let\textit{\ }$\gamma \in \pi _{1}(\Lambda
(n),\ell _{0})$. In view of (\ref{mabeta}) we have the equality 
\begin{equation*}
\mu _{\Lambda }(\gamma ,\ell )=\mu _{\infty }(\gamma \ell _{0,\infty },\ell
_{\infty })-\mu _{\infty }(\ell _{0,\infty },\ell _{\infty })=2m(\gamma )%
\text{.}
\end{equation*}
\end{enumerate}
\end{proof}

The following result allows us to compare the indices $\mu _{\Lambda
}(\lambda ,\ell )$ and $\mu _{\Lambda }(\lambda ,\ell ^{\prime })$
corresponding to different choices of Lagrangians $\ell $ and $\ell ^{\prime
}$.

\begin{proposition}
\label{chgt}For all $\lambda _{12}\in \mathcal{C}(\Lambda (n))$ and $(\ell
,\ell ^{\prime })\in \Lambda (n)^{2}$ we have 
\begin{equation}
\mu _{\Lambda }(\lambda _{12},\ell )-\mu _{\Lambda }(\lambda _{12},\ell
^{\prime })=\func{sign}(\ell _{2},\ell ,\ell ^{\prime })-\func{sign}(\ell
_{1},\ell ,\ell ^{\prime })  \label{diff}
\end{equation}
\end{proposition}

\begin{proof}
Considering the notations of theorem \ref{deep} \ we have 
\begin{multline*}
\mu _{\Lambda }(\lambda _{12},\ell )-\mu _{\Lambda }(\lambda _{12},\ell
^{\prime })= \\
\mu _{\infty }(\ell _{2,\infty },\ell _{\infty })-\mu _{\infty }(\ell
_{2,\infty },\ell _{\infty }^{\prime })-(\mu _{\infty }(\ell _{1,\infty
},\ell _{\infty })-\mu _{\infty }(\ell _{1,\infty },\ell _{\infty }^{\prime
})
\end{multline*}
that is, in view of (\ref{unideux}) 
\begin{multline*}
\mu _{\Lambda }(\lambda _{12},\ell )-\mu _{\Lambda }(\lambda _{12},\ell
^{\prime })= \\
-\mu _{\infty }(\ell _{\infty },\ell _{\infty }^{\prime })+\func{sign}(\ell
_{2},\ell ,\ell ^{\prime })-(-\mu _{\infty }(\ell _{\infty },\ell _{\infty
}^{\prime })+\func{sign}(\ell _{1},\ell ,\ell ^{\prime }))
\end{multline*}
from which (\ref{diff}) immediately follows.
\end{proof}

\subsection{Construction of the canonical index $\bar{\protect\mu}_{\infty }$%
\label{Kashiwa}}

A consequence of the relative uniqueness theorem \ref{un} is that all
Lagrangian intersection indices will be known from the moment we have
determined one of them, and it follows from theorem \ref{deep} that it
suffices for that to construct a Leray index.

Using the intersection theory of Lefschetz chains, Leray has constructed in 
\cite{Leray} (Ch.I, \S 2.5) a function 
\begin{equation*}
m:\Lambda (n)_{\infty }^{2}(0)=\left\{ (\ell _{1,\infty },\ell _{2,\infty
}):\ell _{1}\cap \ell _{2}=0\right\} \longrightarrow \mathbb{Z}
\end{equation*}%
such that 
\begin{equation}
m(\ell _{1,\infty },\ell _{2,\infty })-m(\ell _{1,\infty },\ell _{3,\infty
})+m(\ell _{2,\infty },\ell _{3,\infty })=\func{Inert}(\ell _{1},\ell
_{2},\ell _{3})  \label{leray}
\end{equation}%
where $\func{Inert}(\ell _{1},\ell _{2},\ell _{3})$ is the index of inertia\
of the Lagrangian triple $(\ell _{1},\ell _{2},\ell _{3})$. This index of
inertia is defined in the following way. The transversality condition 
\begin{equation*}
\ell _{1}\cap \ell _{2}=\ell _{2}\cap \ell _{3}=\ell _{3}\cap \ell _{1}=0
\end{equation*}%
being equivalent to 
\begin{equation*}
\ell _{1}\oplus \ell _{2}=\ell _{2}\oplus \ell _{3}=\ell _{3}\oplus \ell
_{1}=X\times X^{\ast }
\end{equation*}%
the relation $z_{1}+z_{2}+z_{3}=0$ ($z_{1}\in \ell _{1}$, $z_{2}\in \ell
_{2} $, $z_{3}\in \ell _{3}$) defines three quadratic forms 
\begin{equation}
z_{1}\longmapsto \omega (z_{2},z_{3})\text{ , }z_{2}\longmapsto \omega
(z_{3},z_{1})\text{ , }z_{3}\longmapsto \omega (z_{1},z_{2})  \label{leray2}
\end{equation}%
such that 
\begin{equation*}
\omega (z_{2},z_{3})=\omega (z_{3},z_{1})=\omega (z_{1},z_{2})\text{.}
\end{equation*}%
These quadratic forms have the same index of inertia and this index is
precisely $\func{Inert}(\ell _{1},\ell _{2},\ell _{3})$.

The function $m$ (which Leray calls \textquotedblleft Maslov
index\textquotedblright ) can be computed using the following formula, due
to Souriau \cite{Souriau2}: 
\begin{equation}
m(\ell _{1,\infty },\ell _{2,\infty })=\frac{1}{2\pi }\left[ \theta
_{1}-\theta _{2}+i\func{Tr}\func{Log}(-w_{1}(w_{2}{}^{-1}))\right] +\frac{n}{%
2}  \label{Souriau}
\end{equation}%
where $\ell _{1,\infty }$ and $\ell _{2,\infty }$ are identified with $%
(w_{1},\theta _{1})$ and $(w_{2},\theta _{2})$, respectively (see\textit{\ }%
\ref{km}). Note that 
\begin{equation}
\dim (\ell _{1}\cap \ell _{2})=\func{corang}(w_{1}-w_{2})  \label{corang1}
\end{equation}%
(it is the multiplicity of the eigenvalue $+1$ of $w_{1}(w_{2}^{-1})$). The
logarithm (\ref{Souriau}) is defined , for every matrix $u$ that has not the
eigenvalue $-1$, by the formula 
\begin{equation*}
\limfunc{Log}(u)=\int_{-\infty }^{0}\left( \left( \lambda I-u\right)
^{-1}-\left( \lambda -1\right) ^{-1}I\right) \,d\lambda
\end{equation*}%
and hence $\func{Log}w_{1}(w_{2}^{-1})$ is well-defined since in view of (%
\ref{corang1}) 
\begin{equation*}
\ell _{1}\cap \ell _{2}=0\Longleftrightarrow \det (w_{1}-w_{2})\neq 0,
\end{equation*}%
that is 
\begin{equation*}
+1\text{\ \textit{is not an eigenvalue of} }w_{1}(w_{2}^{-1})\text{.}
\end{equation*}%
The function $m$ possesses the following properties:

\begin{proposition}
(1) $m$ is locally constant on its domain 
\begin{equation*}
\Lambda _{\infty }^{2}(n;0)=\left\{ (\ell _{1,\infty },\ell _{2,\infty
}):\ell _{1}\cap \ell _{2}=0\right\} \text{;}
\end{equation*}
(2) The action of $\pi _{1}(\Lambda (n))$ on $m$ is given by 
\begin{equation}
m(\gamma _{1}\ell _{1,\infty },\gamma _{2}\ell _{2,\infty })=m(\ell
_{1,\infty },\ell _{2,\infty })+m(\gamma _{1})-m(\gamma _{2})\text{.}
\label{masbeta}
\end{equation}
\end{proposition}

\begin{proof}
(1): Souriaus' formula (\ref{Souriau}) shows that $m$ is continuous and
hence locally constant (it takes its values in a discrete space). (2): To
prove the equality (\ref{masbeta}) we first note that $\pi _{1}(\Lambda
(n))\equiv \pi _{1}(W(n))$ can be identified with the group 
\begin{equation*}
G=\left\{ (I,2k\pi ):k\in \mathbb{Z}\right\}
\end{equation*}%
acting on $W_{\infty }(n)$ by 
\begin{equation*}
(I,2k\pi )(w,\theta )=(w,\theta +2k\pi )\text{.}
\end{equation*}%
Under this identification, the generator $\beta $ of $\pi _{1}(\Lambda (n))$
defined in Proposition \ref{un} is identified with $(I,2\pi )$. For all $%
\gamma \in \pi _{1}(\Lambda (n))$ there exists an integer $k$ such that $%
\gamma =\beta ^{k}$, and hence 
\begin{equation*}
\gamma \ell _{\infty }=\beta ^{k}\ell _{\infty }=(w,\theta +2k\pi )\text{ \ 
\textit{if} \ }\ell _{\infty }=(w,\theta )\text{.}
\end{equation*}%
In view of the definition (\ref{Souriau}) of $m$ we thus obtain 
\begin{align*}
m(\gamma _{1}\ell _{1,\infty },\gamma _{2}\ell _{2,\infty })& =m(\beta
^{k_{1}}\ell _{1,\infty },\beta ^{k_{2}}\ell _{2,\infty })-m(\ell _{1,\infty
},\ell _{2,\infty }) \\
& =\frac{1}{2\pi }\left[ 2k_{1}\pi -2k_{2}\pi \right] \\
& =k_{1}-k_{2}
\end{align*}%
hence (\ref{masbeta}) in view of the formulae (\ref{mabetaun}) and (\ref%
{maprop}).
\end{proof}

Let us now set 
\begin{equation}
\tau =2\func{Inert}-n\text{.}  \label{taui}
\end{equation}%
We have the important relation 
\begin{equation}
\tau (\ell _{1},\ell _{2},\ell _{3})=\tau ^{_{+}}-\tau ^{_{-}}  \label{taua}
\end{equation}%
where $\tau ^{_{+}}$ (resp. $\tau ^{_{-}}$) is the number of $>0$ (resp. $<0$%
) eigenvalues of the quadratic form 
\begin{equation}
Q(z_{1},z_{2},z_{3})=\omega (z_{1},z_{2})+\omega (z_{2},z_{3})+\omega
(z_{3},z_{1})  \label{kubai}
\end{equation}%
(see (see \cite{JMPA}). It follows that $\tau $ is the signature of
Demazure--Kashiwara (see \cite{Demazure,Marle,LV}). It is defined on all
triples $(\ell _{1},\ell _{2},\ell _{3}),$ which is not the case for $\func{%
Inert}$ as seen above.

\begin{remark}
The function $\tau $ is often called ``Maslov index'', or ``trilateral
Maslov index'' in the literature (for instance in \cite{CLM} and \cite{LV}).
We will not use this terminology.
\end{remark}

The signature of Demazure--Kashiwara has the following properties (see
Libermann-Marle \cite{Marle} or Lion-Vergne \cite{LV}; in Capell \textit{et
al}. \cite{CLM} one can find other interesting properties).

\begin{proposition}
We have:

\begin{enumerate}
\item[(1)] \label{tautau} $\tau \in C_{ant}^{2}(\Lambda (n),\mathbb{Z})^{Sp}$%
, i.e. , $\tau $ is a $2$-cocycle, totally antisymmetric and $Sp(n)$%
-invariant : \ \ 
\begin{equation}
\partial \tau =0\text{ \ , }\varepsilon ^{\ast }\tau =(-1)^{sgn(\varepsilon
)}\tau \text{\ \ , \ }s^{\ast }\tau =\tau  \label{tauco}
\end{equation}%
($\varepsilon $ all permutation of $\left\{ 1,2,3\right\} $, $s\in Sp$).

\item[(2)] $\tau $ is locally constant on the sets 
\begin{equation*}
\left\{ (\ell _{1},\ell _{2},\ell _{3}):\dim (\ell _{1}\cap \ell
_{2})=k_{1},\dim (\ell _{1}\cap \ell _{2})=k_{2},\dim (\ell _{1}\cap \ell
_{2})=k_{3}\right\}
\end{equation*}

\item[(3)] Let $A$ a symmetric matrix of order $n$ and $\ell _{A}$ the
Lagrangian with equation $p=Ax$. We have 
\begin{equation}
\tau (X^{\ast },\ell _{A},X)=\func{sign}A  \label{taula}
\end{equation}
where $\func{sign}A$ is the difference between the number of eigenvalues $>0$
and $<0$ of $A$.\textit{\medskip }
\end{enumerate}
\end{proposition}

Let us now set 
\begin{equation}
\bar{\mu}_{\infty }(\ell _{1,\infty },\ell _{2,\infty })=2m(\ell _{1,\infty
},\ell _{2,\infty })-n  \label{zet}
\end{equation}
that is, in view of (\ref{Souriau}) : 
\begin{equation}
\bar{\mu}_{\infty }(\ell _{1,\infty },\ell _{2,\infty })=\frac{1}{\pi }\left[
\theta _{1}-\theta _{2}+i\func{Tr}\func{Log}(-w_{1}(w_{2}^{-1}))\right] 
\text{.}  \label{souribis}
\end{equation}
For $(\ell _{1,\infty },\ell _{2,\infty })\in \Lambda _{\infty }(n)^{2}$ let
us choose $\ell _{3,\infty }\in \Lambda (n)_{\infty }$ such that 
\begin{equation}
\ell _{1}\cap \ell _{3}=\ell _{2}\cap \ell _{3}=0  \label{123}
\end{equation}
and \textit{define} 
\begin{equation}
\bar{\mu}_{\infty }(\ell _{1,\infty },\ell _{2,\infty })=\bar{\mu}_{\infty
}(\ell _{1,\infty },\ell _{3,\infty })-\bar{\mu}_{\infty }(\ell _{2,\infty
},\ell _{3,\infty })+\tau (\ell _{1},\ell _{2},\ell _{3})\text{.}
\label{mules}
\end{equation}

We can verify, using the cocycle property of $\ \tau $ that

(1) The right-hand side of (\ref{mules}) does not depend on the choice of $%
\ell _{3,\infty }$ such that we have (\ref{123}), which justifies the
notation $\bar{\mu}_{\infty }(\ell _{1,\infty },\ell _{2,\infty })$ (see 
\cite{JMPA} where $\bar{\mu}_{\infty }$ is denoted by $\mu $).

(2) We have $\partial \bar{\mu}_{\infty }=\pi ^{\ast }\tau $, i.e., 
\begin{equation}
\bar{\mu}_{\infty }(\ell _{1,\infty },\ell _{2,\infty })-\bar{\mu}_{\infty
}(\ell _{1,\infty },\ell _{3,\infty })+\bar{\mu}_{\infty }(\ell _{2,\infty
},\ell _{3,\infty })=\tau (\ell _{1},\ell _{2},\ell _{3})  \label{remule}
\end{equation}
for all $(\ell _{1,\infty },\ell _{2,\infty },\ell _{3,\infty })\in \Lambda
_{\infty }(n)^{3}$.\medskip

\begin{definition}
We will call $\bar{\mu}_{\infty }$ the ``canonical Leray index'' and the
Lagrangian and symplectic intersection indices associated to $\bar{\mu}%
_{\infty }$\ will be call the ``canonical intersection indices''. These will
be denoted $\bar{\mu}_{\Lambda }$ and $\bar{\mu}_{Sp}$, respectively.
\end{definition}

The following theorem will show that $\bar{\mu}_{\infty }$ is indeed a Leray
index.

\begin{theorem}
\label{fonda} The function $\bar{\mu}_{\infty }$ defined by (\ref{mules}) is
the unique $1$-cochain on $\Lambda _{\infty }(n)$ having the following two
properties:

\begin{enumerate}
\item[(1)] $\partial \bar{\mu}_{\infty }=\pi ^{\ast }\tau $.

\item[(2)] $\mu _{\infty }$ is locally constant on the sets $\Lambda
_{\infty }^{2}(n;k)$.
\end{enumerate}

Furthermore, this $1$-cochain has the following property:

\begin{enumerate}
\item[(3)] The action of $\pi _{1}(\Lambda (n))$ on $\bar{\mu}_{\infty }$ is
given by 
\begin{equation}
\bar{\mu}_{\infty }(\gamma _{1}\ell _{1,\infty },\gamma _{2}\ell _{2,\infty
})=\bar{\mu}_{\infty }(\ell _{1,\infty },\ell _{2,\infty })+2(m(\gamma
_{1})-m(\gamma _{2}))\text{.}  \label{entiers}
\end{equation}
\end{enumerate}
\end{theorem}

\begin{proof}
Let us begin by showing that there exists at most one $1$-cochain $\mu
_{\infty }$ having the properties (1) and (2). To do so, let us denote by $%
\nu $ the difference $\bar{\mu}_{\infty }-\bar{\mu}_{\infty }^{\prime }$
between two cochains having the stated properties. In this case 
\begin{equation*}
\nu (\ell _{1,\infty },\ell _{2,\infty })=\nu (\ell _{1,\infty },\ell
_{3,\infty })-\nu (\ell _{2,\infty },\ell _{3,\infty })\text{.}
\end{equation*}
for all $(\ell _{1,\infty },\ell _{2,\infty },\ell _{3,\infty })$. We choose 
$\ell _{3}$ in such a way that $\ell _{1}\cap \ell _{3}=\ell _{2}\cap \ell
_{3}=0$. In view of the definition of $\nu $ and the axiom (L$_{3}$), $\nu
(\ell _{1,\infty },\ell _{2,\infty })$ will remain constant in connected
neighborhoods of $\ell _{1,\infty }$ and $\ell _{2,\infty }$. The function $%
\nu $ is thereby locally constant on $\Lambda _{\infty }(n)^{2}$. But this
space is connected, and having realized this we conclude that $\nu $ is in
fact constant. Also choosing $\ell _{1,\infty }=\ell _{2,\infty }$ we find
that the values of $\nu $ is zero, and hence that $\bar{\mu}_{\infty }=\bar{%
\mu}_{\infty }^{\prime }$.

That $\mu _{\infty }$ possesses property (1) is clear from its construction.
Let us therefore concentrate on the proof of property (2). Suppose first
that $k=0$. In view of the Souriau formula (\ref{Souriau}) and definition (%
\ref{zet}), the function $\bar{\mu}_{\infty }$ is continuous on $\Lambda
_{\infty }^{2}(k)$ and is thereby locally constant since $\bar{\mu}_{\infty
} $ takes discrete values. For arbitrary $k$, it is enough to choose $\ell
_{3} $ transversal to $\ell _{1}$ and $\ell _{2}$ in (\ref{mules}) in such a
way that $\bar{\mu}_{\infty }(\ell _{1,\infty },\ell _{3,\infty })$ and $%
\bar{\mu}_{\infty }(\ell _{2,\infty },\ell _{3,\infty })$ remain constant,
and to apply the property (2) of $\tau $ in Proposition \ref{tautau}. Let us
finally show that $\mu _{\infty }$ has property (3). Choosing once more $%
\ell _{3}$ such that $\ell _{1}\cap \ell _{3}=\ell _{2}\cap \ell _{3}=0$, we
find that in view of (\ref{mules}) and the fact that for all $\gamma \in \pi
_{1}(\Lambda (n))$ the element $\gamma \ell _{\infty }$ has the same
projection $\ell $ as $\ell _{\infty }$ 
\begin{multline*}
\bar{\mu}_{\infty }(\gamma _{1}\ell _{1,\infty },\gamma _{2}\ell _{2,\infty
})=\bar{\mu}_{\infty }(\gamma _{1}\ell _{1,\infty },\ell _{3,\infty }) \\
-\bar{\mu}_{\infty }(\gamma _{2}\ell _{2,\infty },\ell _{3,\infty })+\tau
(\ell _{1},\ell _{2},\ell _{3})\text{.}
\end{multline*}
By (\ref{masbeta}) and the definition (\ref{zet}) of $\bar{\mu}_{\infty }$
in the transversal case, we obtain 
\begin{equation*}
\left\{ 
\begin{array}{c}
\bar{\mu}_{\infty }(\gamma _{1}\ell _{1,\infty },\ell _{3,\infty })=\bar{\mu}%
_{\infty }(\ell _{1,\infty },\ell _{2,\infty })+2m(\gamma _{1})\smallskip \\ 
\bar{\mu}_{\infty }(\gamma _{2}\ell _{2,\infty },\ell _{3,\infty })=\bar{\mu}%
_{\infty }(\ell _{2,\infty },\ell _{3,\infty })+2m(\gamma _{2})%
\end{array}
\right.
\end{equation*}
and (\ref{entiers}) follows.
\end{proof}

We point out the fact that Theorem \ref{fonda} also shows that:

\begin{corollary}
The Demazure--Kashiwara index is identical to the signature cocycle, i.e., $%
\tau =\func{sign}$.
\end{corollary}

\begin{proof}
This is clear since $\mu _{\infty }$ is a Leray index, so $\partial \mu
_{\infty }=\func{sign}$, and since we also have $\partial \mu _{\infty
}=\tau $.
\end{proof}

\section{Further properties of the intersection indices}

We shall see in the following paragraphs that the properties of the
canonical Leray index enable us to show several properties of the Lagrangian
and symplectic intersection indices. Indeed, we shall for instance see that
these intersection indices are symplectic invariant.

\subsection{ symplectic invariance}

As we mentioned the results of the preceding section show that the
Lagrangian and symplectic intersection indices $\mu _{\Lambda (n)}$ and $\mu
_{Sp}$, respectively, are symplectic invariants. In order to prove this fact
we need first to show that the Leray index is invariant under the action of
the universal covering space of the symplectic group the following

\begin{lemma}
\label{cucul}The Leray index is invariant under the action of the universal
covering space $Sp_{\infty }(n)$ of $Sp(n)$: 
\begin{equation}
\mu _{\infty }(s_{\infty }\ell _{1,\infty },s_{\infty }\ell _{2,\infty
})=\mu _{\infty }(\ell _{1,\infty },\ell _{2,\infty })  \label{leraysp}
\end{equation}
for all $(s_{\infty },\ell _{1,\infty },\ell _{2,\infty })\in Sp_{\infty
}(n)\times \Lambda (n)_{\infty }\times \Lambda (n)_{\infty }$.
\end{lemma}

\begin{proof}
Consider the mapping 
\begin{equation*}
\mu _{\infty }^{\prime }:(\ell _{1,\infty },\ell _{2,\infty })\longmapsto
\mu _{\infty }(s_{\infty }\ell _{1,\infty },s_{\infty }\ell _{2,\infty })%
\text{.}
\end{equation*}
Since $\tau $ is invariant under the action of $Sp(n)$ we have $\partial \mu
_{\infty }^{\prime }=\tau $. Of course, the function $\mu _{\infty }^{\prime
}$ is locally constant on $\Lambda _{\infty }^{2}(n;0)$ and is therefore
identical to $\mu _{\infty }$ in view of the uniqueness property (1) in
Theorem \ref{fonda}.
\end{proof}

Having shown this we may now state the following essential result. The
invariance of the symplectic intersection indices will then follow from
their definition.

\begin{theorem}
The Lagrangian intersection indices are symplectic invariants, i.e., 
\begin{equation}
\mu _{\Lambda }(s\lambda _{12},s\ell )=\mu _{\Lambda }(\lambda _{12},\ell )%
\text{ \ for all \ }s\in Sp(n)\text{.}  \label{invsymp}
\end{equation}
\end{theorem}

\begin{proof}
In view of Theorem \ref{un} and using the fact that 
\begin{equation*}
\dim (s\ell _{1}\cap s\ell _{2})=\dim (\ell _{1}\cap \ell _{2})
\end{equation*}
for all $(s,\ell _{1},\ell _{2})\in Sp(n)\times \Lambda (n)^{2}$ it suffices
to consider the case in which $\mu _{\Lambda }$ is the canonical
intersection index $\bar{\mu}_{\Lambda }$. for $s\in Sp(n)$ we consider an
arbitrary symplectic path $\sigma \in \mathcal{C}(Sp(n))$ joining the
identity $I$ of $Sp(n)$ to $s$. We denote $s_{\infty }$ the homotopy
equivalence class of that path, $s_{\infty }\in Sp_{\infty }(n)$. As in
Theorem \ref{deep} we denote by $\lambda _{01}$ a path in $\Lambda (n)$
joining $\ell _{0}$ to $\ell _{1}$ and $\ell _{1,\infty }$ it homotopy
class. $\ell _{2,\infty }$ is here the homotopy class of the path $\lambda
_{01}\ast \lambda _{12}$. Denoting by $\lambda $ an arbitrary path joining $%
\ell _{0}$ to $\ell $ in $\Lambda (n)$, and $\ell _{\infty }$ its homotopy
class we find 
\begin{equation*}
s_{\infty }\ell _{1,\infty }=\text{class}\left[ t\longmapsto \sigma
(t)\lambda _{01}(t),0\leq t\leq 1\right]
\end{equation*}
\begin{equation*}
s_{\infty }\ell _{2,\infty }=\text{class}\left[ t\longmapsto \left\{ 
\begin{array}{c}
\sigma (t)\lambda _{01}(t),0\leq t\leq \tfrac{1}{2}\smallskip \\ 
\sigma (2t-1)\lambda _{12}(2t-1),0\leq t\leq \tfrac{1}{2}%
\end{array}
\right. \right]
\end{equation*}
\begin{equation*}
s_{\infty }\ell _{\infty }=\text{class}\left[ t\longmapsto \sigma (t)\lambda
(t),0\leq t\leq 1\right] \text{.}
\end{equation*}
Thus, by (\ref{foufoun}) in Theorem \ref{deep} we obtain 
\begin{equation*}
\bar{\mu}_{\infty }(s_{\infty }\ell _{2,\infty },s_{\infty }\ell _{\infty })-%
\bar{\mu}_{\infty }(s_{\infty }\ell _{1,\infty },s_{\infty }\ell _{\infty })=%
\bar{\mu}_{\Lambda }(s\lambda _{01},s\ell )\text{.}
\end{equation*}
But the left-hand side of this equality is 
\begin{equation*}
\bar{\mu}_{\infty }(\ell _{2,\infty },\ell _{\infty })-\bar{\mu}_{\infty
}(\ell _{1,\infty },\ell _{\infty })=\bar{\mu}_{\Lambda }(\lambda _{01},\ell
)
\end{equation*}
by Lemma \ref{cucul}, which proves (\ref{invsymp}).
\end{proof}

\subsection{The indices $\bar{\protect\mu}_{\ell }$ on $Sp_{\infty }$}

Theorem \ref{deep} expresses every Lagrangian intersection index as the
difference between two values of the Leray index on $\Lambda _{\infty }(n)$.
We will now show that there is a similar result for symplectic intersection
indices. We begin by showing the following result (\textit{cf}. \cite{JMPA}):

\begin{lemma}
Let $s_{\infty }\in Sp_{\infty }(n)$ and $\ell _{\infty }\in \Lambda
_{\infty }(n)$, with projection $\ell \in \Lambda (n)$. The integer $\bar{\mu%
}_{\infty }(s_{\infty }\ell _{\infty },\ell _{\infty })$ depends solely on $%
(s_{\infty },\ell )$, and not on the choice of the projection $\ell $ of $%
\ell _{\infty }$. Hence, there exists for every $\ell \in \Lambda (n)$, a
function $\bar{\mu}_{\ell }:Sp_{\infty }(n)\longrightarrow \mathbb{Z}$ \
such that 
\begin{equation*}
\bar{\mu}_{\ell }(s_{\infty })=\bar{\mu}_{\infty }(s_{\infty }\ell _{\infty
},\ell _{\infty })\text{.}
\end{equation*}
\end{lemma}

\begin{proof}
\smallskip Suppose that $\ell _{\infty }$ and $\ell _{\infty }^{\prime }$
are two elements of $\Lambda _{\infty }(n)$ having the same projection $\ell
\in \Lambda (n)$. There is then a $k\in \mathbb{Z}$ such that $\ell _{\infty
}^{\prime }=\beta ^{k}\ell _{\infty }$. ($\beta $ being the generator of $%
\pi _{1}(\Lambda (n))$ whose natural image in $\mathbb{Z}$ is $+1$), and in
view of (\ref{alfabeta}) we obtain 
\begin{eqnarray*}
\bar{\mu}_{\infty }(s_{\infty }\ell _{\infty }^{\prime },\ell _{\infty
}^{\prime })=\bar{\mu}_{\infty }(s_{\infty }(\beta ^{k}\ell _{\infty
}),\beta ^{k}\ell _{\infty }) \\
=\bar{\mu}_{\infty }(\beta ^{k}(s_{\infty }\ell _{\infty }),\beta ^{k}\ell
_{\infty })
\end{eqnarray*}
that is 
\begin{equation*}
\bar{\mu}_{\infty }(s_{\infty }\ell _{\infty }^{\prime },\ell _{\infty
}^{\prime })=\bar{\mu}_{\infty }(s_{\infty }\ell _{\infty },\ell _{\infty })
\end{equation*}
by the property (\ref{entiers}) in Theorem \ref{fonda}) of $\bar{\mu}%
_{\infty }$.
\end{proof}

We will call the index $\bar{\mu}_{\ell }$ the ``Leray index on $Sp_{\infty
}(n)$'' relative to the Lagrangian $\ell $. We will consider in what follows 
$s_{\infty }\in Sp_{\infty }(n)$ as the homotopy class of a continuous path
joining the identity $I$ to an element $s$ (the projection of $s_{\infty }$)
in $Sp(n)$.

\begin{proposition}
Let $\sigma _{12}\in \mathcal{C}(Sp(n))$ be a symplectic path joining $s_{1}$
to $s_{2}$ in $Sp(n)$. Let $s_{1,\infty }$ be an arbitrary element of $%
Sp_{\infty }(n)$ covering $s_{1}$ and $s_{2,\infty }$ the homotopy class of $%
\sigma _{01}\ast \sigma _{12}$ ($\sigma _{01}$ a representative of $%
s_{1,\infty }$). We have 
\begin{equation}
\bar{\mu}_{Sp}(\sigma _{12},\ell )=\bar{\mu}_{\ell }(s_{2,\infty })-\bar{\mu}%
_{\ell }(s_{1,\infty })\text{.}  \label{dire}
\end{equation}
\end{proposition}

\begin{proof}
Consider $\ell _{\infty }$ to be the homotopy class of an arbitrary path $%
\lambda $ joining $\ell _{0}$ (the base point of $\Lambda (n)_{\infty }$) to 
$\ell $. We have 
\begin{equation*}
s_{1,\infty }\ell _{\infty }=\text{class}\left[ t\longmapsto \sigma
_{01}(t)\lambda (t),0\leq t\leq 1\right]
\end{equation*}
and 
\begin{equation*}
s_{2,\infty }\ell _{\infty }=\text{class}\left[ 
\begin{array}{c}
t\longmapsto \sigma _{01}(2t)\lambda (2t),0\leq t\leq \frac{1}{2} \\ 
t\longmapsto \sigma _{12}(2t-1)\lambda (2t-1),\frac{1}{2}\leq t\leq 1%
\end{array}
\right]
\end{equation*}
hence, by the formulae (\ref{formuleun}) and (\ref{foufoun}) one obtains 
\begin{align*}
\bar{\mu}_{Sp}(\sigma _{12},\ell )=\bar{\mu}_{\Lambda (n)}(\sigma _{12}\ell
,\ell ) \\
=\bar{\mu}_{\infty }(s_{2,\infty }\ell _{\infty },\ell _{\infty })-\bar{\mu}%
_{\infty }(s_{1,\infty }\ell _{\infty },\ell _{\infty })
\end{align*}
that is (\ref{dire}).
\end{proof}

The two important properties of $\bar{\mu}_{\ell ,\infty }$ below are
consequences of cohomological property (\ref{remule}) of the Leray index

\begin{proposition}
We have

\begin{enumerate}
\item[(1)] 
\begin{equation}
\bar{\mu}_{\ell }(s_{1,\infty }s_{2,\infty })=\bar{\mu}_{\ell }(s_{1,\infty
})+\bar{\mu}_{\ell }(s_{2,\infty })+\tau (\ell ,s_{1}\ell ,s_{1}s_{2}\ell )
\label{cha}
\end{equation}
for all $(s_{1,\infty },s_{2,\infty })\in Sp_{\infty }(n)$.

\item[(2)] Let $\ell $ and $\ell ^{\prime }$ be two Lagrangian, then 
\begin{equation}
\bar{\mu}_{\ell }(s_{\infty })-\bar{\mu}_{\ell ^{\prime }}(s_{\infty })=\tau
(s\ell ,\ell ,\ell ^{\prime })-\tau (s\ell ,s\ell ^{\prime },\ell ^{\prime })%
\text{.}  \label{sl}
\end{equation}
\end{enumerate}
\end{proposition}

\begin{proof}
(1) By the definition of $\bar{\mu}_{\ell }$ we see that 
\begin{multline*}
\bar{\mu}_{\ell }(s_{1,\infty }s_{2,\infty })-\bar{\mu}_{\ell }(s_{1,\infty
})-\bar{\mu}_{\ell }(s_{2,\infty })= \\
\bar{\mu}_{\infty }(s_{1,\infty }s_{2,\infty }\ell _{\infty },\ell _{\infty
})-\bar{\mu}_{\infty }(s_{1,\infty }\ell _{\infty },\ell _{\infty })-\bar{\mu%
}_{\infty }(s_{2,\infty }\ell _{\infty },\ell _{\infty })
\end{multline*}
that is, using the $Sp_{\infty }(n)$--invariance and the antisymmetry of $%
\bar{\mu}_{\infty }$: 
\begin{multline*}
\bar{\mu}_{\ell }(s_{1,\infty }s_{2,\infty })-\bar{\mu}_{\ell }(s_{1,\infty
})-\bar{\mu}_{\ell }(s_{2,\infty })=\bar{\mu}_{\infty }(s_{1,\infty
}s_{2,\infty }\ell _{\infty },\ell _{\infty })+ \\
\bar{\mu}_{\infty }(\ell _{\infty },s_{1,\infty }\ell _{\infty })-\bar{\mu}%
_{\infty }(s_{1,\infty }s_{2,\infty }\ell _{\infty },s_{1,\infty }\ell
_{\infty })\text{.}
\end{multline*}
In view of the property (\ref{remule}) of $\bar{\mu}_{\infty }$ we obtain 
\begin{align*}
\bar{\mu}_{\ell }(s_{1,\infty }s_{2,\infty })-\bar{\mu}_{\ell }(s_{1,\infty
})-\bar{\mu}_{\ell }(s_{2,\infty })=\tau (s_{1}s_{2}\ell ,\ell ,s_{1}\ell )
\\
=\tau (\ell ,s_{1}\ell ,s_{1}s_{2}\ell )\text{.}
\end{align*}
(2) By (\ref{remule}) and the $Sp_{\infty }(n)$--invariance of $\bar{\mu}%
_{\infty }$, 
\begin{equation*}
\bar{\mu}_{\infty }(s_{\infty }\ell _{\infty },\ell _{\infty })-\bar{\mu}%
_{\infty }(s_{\infty }\ell _{\infty },\ell _{\infty }^{\prime })+\bar{\mu}%
_{\infty }(s_{\infty }\ell _{\infty },s_{\infty }\ell _{\infty }^{\prime
})=\tau (s\ell ,\ell ,\ell ^{\prime })
\end{equation*}
as well as 
\begin{equation*}
\bar{\mu}_{\infty }(s_{\infty }\ell _{\infty },s_{\infty }\ell _{\infty
}^{\prime })-\bar{\mu}_{\infty }(s_{\infty }\ell _{\infty },\ell _{\infty
}^{\prime })+\bar{\mu}_{\infty }(s_{\infty }\ell _{\infty }^{\prime },\ell
_{\infty }^{\prime })=\tau (s\ell ,s\ell ^{\prime },\ell ^{\prime })
\end{equation*}
gives us (\ref{sl}) by substracting the first identity from the second.
\end{proof}

\begin{corollary}
Given two arbitrary Lagrangians $\ell $ and $\ell ^{\prime }$ we have 
\begin{eqnarray}
\bar{\mu}_{Sp}(\sigma _{12},\ell )-\bar{\mu}_{Sp}(\sigma _{12},\ell ^{\prime
}) &=&\tau (s_{2}\ell ,\ell ,\ell ^{\prime })-\tau (s_{2}\ell ,s_{2}\ell
^{\prime },\ell ^{\prime })  \label{hum} \\
&&-(\tau (s_{1}\ell ,\ell ,\ell ^{\prime })-\tau (s_{1}\ell ,s_{1}\ell
^{\prime },\ell ^{\prime }))  \notag
\end{eqnarray}
where $s_{1}$ and $s_{2}$ are the endpoints of the path $\sigma _{12}$.
\end{corollary}

\begin{proof}
Formula (\ref{dire}) yields 
\begin{multline*}
\bar{\mu}_{Sp}(\sigma _{12},\ell )-\bar{\mu}_{Sp}(\sigma _{12},\ell ^{\prime
})= \\
\bar{\mu}_{\ell }(s_{2,\infty })-\bar{\mu}_{\ell }(s_{1,\infty })-(\bar{\mu}%
_{\ell ^{\prime }}(s_{2,\infty })-\bar{\mu}_{\ell ^{\prime }}(s_{1,\infty }))
\end{multline*}
hence (\ref{hum}), using(\ref{sl}).
\end{proof}

\subsection{Dimensional additivity}

The canonical index $\bar{\mu}_{\Lambda }$ has the interesting ``dimensional
additivity'' property, which is often included as an axiom for the
intersection indices (\textit{cf}. \cite{Duff}).

Let us introduces the following notations. Let $n^{\prime }$ and $n^{\prime
\prime }$ be two integers $>0$ and set $n=n^{\prime }+n^{\prime \prime }$.
We denote by $\Lambda (n^{\prime })$ and $\Lambda (n^{\prime \prime })$
(resp. $Sp(n^{\prime })$ and $Sp(n^{\prime \prime })$) the Lagrangian
Grassmannians (resp. the symplectic groups) corresponding to the symplectic
spaces $(X^{\prime }\times X^{\prime \ast },\omega ^{\prime })$ and $%
(X^{\prime \prime }\times X^{\prime \prime \ast },\omega ^{\prime \prime })$%
. The direct sum 
\begin{equation*}
\Lambda (n^{\prime })\oplus \Lambda (n^{\prime \prime })=\left\{ \ell
^{\prime }\oplus \ell ^{\prime \prime }:\ell ^{\prime }\in \Lambda
(n^{\prime }),\ell ^{\prime \prime }\in \Lambda (n^{\prime \prime })\right\}
\end{equation*}%
is identified with a submanifold of $\Lambda (n)$ and 
\begin{equation*}
Sp(n^{\prime })\oplus Sp(n^{\prime \prime })=\left\{ s^{\prime }\oplus
s^{\prime \prime }:s^{\prime }\in Sp(n^{\prime }),s^{\prime \prime }\in
Sp(n^{\prime \prime })\right\}
\end{equation*}%
with a subgroup of $Sp(n)$. By definition 
\begin{equation}
(s^{\prime }\oplus s^{\prime \prime })(z^{\prime }\oplus z^{\prime \prime
})=s^{\prime }(z^{\prime })\oplus s^{\prime \prime }(z^{\prime \prime })%
\text{.}  \label{sprime}
\end{equation}%
Denoting $\tau ^{\prime }$ and $\tau ^{\prime \prime }$ the signatures of
triples of elements of $\Lambda (n^{\prime })$ and $\Lambda (n^{\prime
\prime })$ we immediately verify that 
\begin{equation}
\tau (\ell _{1}^{\prime }\oplus \ell _{1}^{\prime \prime },\ell _{2}^{\prime
}\oplus \ell _{2}^{\prime \prime },\ell _{3}^{\prime }\oplus \ell
_{3}^{\prime \prime })=\tau ^{\prime }(\ell _{1}^{\prime },\ell _{2}^{\prime
},\ell _{3}^{\prime })+\tau ^{\prime \prime }(\ell _{1}^{\prime \prime
},\ell _{2}^{\prime \prime },\ell _{3}^{\prime \prime })\text{.}
\label{tausomme}
\end{equation}

In view of the identification of $\Lambda _{\infty }(n)$ with $W_{\infty
}(n) $ we may define 
\begin{equation*}
\ell _{\infty }^{\prime }\oplus \ell _{\infty }^{\prime \prime }=\left(
w^{\prime }\oplus w^{\prime \prime },\theta ^{\prime }+\theta ^{\prime
\prime }\right)
\end{equation*}
if $\ell _{\infty }^{\prime }=\left( w^{\prime },\theta ^{\prime }\right) $
and $\ell _{\infty }^{\prime \prime }=\left( w^{\prime \prime },\theta
^{\prime \prime }\right) $. This defines an element of $\Lambda (n)_{\infty
} $ since $\left( w^{\prime }\oplus w^{\prime \prime },\theta ^{\prime
}+\theta ^{\prime \prime }\right) \in W_{\infty }(n)$.

The following result describes the relation between the Leray indices on the
Maslov bundles $\Lambda _{\infty }(n^{\prime })$, $\Lambda _{\infty
}(n^{\prime \prime })$ and the Leray index on $\Lambda _{\infty }(n)$

\begin{proposition}
Let $\bar{\mu}_{\infty }^{\prime }$ and $\bar{\mu}_{\infty }^{\prime \prime
} $ be the Leray indices on $\Lambda _{\infty }(n^{\prime })$ and $\Lambda
_{\infty }(n^{\prime \prime })$, respectively. The relation between the
Leray indices $\bar{\mu}_{\infty }^{\prime }$ and $\bar{\mu}_{\infty
}^{\prime \prime }$ and the Leray index $\bar{\mu}_{\infty }$ on $\Lambda
_{\infty }(n)$ is given by 
\begin{equation}
\bar{\mu}_{\infty }\left( \ell _{1,\infty }^{\prime }\oplus \ell _{1,\infty
}^{\prime \prime },\ell _{2,\infty }^{\prime }\oplus \ell _{2,\infty
}^{\prime \prime }\right) =\bar{\mu}_{\infty }^{\prime }\left( \ell
_{1,\infty }^{\prime },\ell _{2,\infty }^{\prime }\right) +\bar{\mu}_{\infty
}^{\prime \prime }\left( \ell _{1,\infty }^{\prime \prime },\ell _{2,\infty
}^{\prime \prime }\right) \text{.}  \label{muso}
\end{equation}
\end{proposition}

\begin{proof}
Let us begin by supposing that $\ell _{1,\infty }^{\prime },\ell _{2,\infty
}^{\prime }\in \Lambda _{\infty }(n^{\prime })$ and $\ell _{1,\infty
}^{\prime \prime }$, $\ell _{2,\infty }^{\prime \prime }\in \Lambda _{\infty
}(n^{\prime \prime })$ are transversal: $\ell _{1}^{\prime }\cap \ell
_{2}^{\prime }=0$ and $\ell _{1}^{\prime \prime }\cap \ell _{2}^{\prime
\prime }=0$. Then $\ell _{1,\infty }^{\prime }\oplus \ell _{1,\infty
}^{\prime \prime }$ and $\ell _{2,\infty }^{\prime }\oplus \ell _{2,\infty
}^{\prime \prime }$ are also transversal: 
\begin{equation*}
(\ell _{1}^{\prime }\oplus \ell _{1}^{\prime \prime })\cap (\ell
_{2}^{\prime }\oplus \ell _{2}^{\prime \prime })=0
\end{equation*}
and one shows that 
\begin{multline*}
\limfunc{Log}(\left( w_{1}^{\prime }\oplus w_{1}^{\prime \prime }\right)
\left( w_{2}^{\prime }\oplus w_{2}^{\prime \prime }\right) ^{-1})= \\
\limfunc{Log}(w_{1}^{\prime }(w_{2}^{\prime })^{-1})\oplus \limfunc{Log}%
(w_{1}^{\prime \prime }(w_{2}^{\prime \prime })^{-1})
\end{multline*}
under the identifications $\ell _{1}^{\prime }\equiv w_{1}^{\prime }$, $\ell
_{2}^{\prime }\equiv w_{2}^{\prime }$, etc. (see \cite{SdG}). Formula (\ref%
{muso}) follows in the considered case. The general case can be reduced to
the transversal case by using the definition (\ref{mules}) and the formula (%
\ref{tausomme}).
\end{proof}

From this Proposition we immediately deduce that:

\begin{corollary}
Let $\lambda ^{\prime }\in C(\Lambda (n^{\prime }))$ and $\lambda ^{\prime
\prime }\in C(\Lambda (n^{\prime \prime }))$ be Lagrangian paths in $\Lambda
(n^{\prime })$ and $\Lambda (n^{\prime \prime })$. Then 
\begin{equation*}
\lambda ^{\prime }\oplus \lambda ^{\prime \prime }\in \mathcal{C}(\Lambda
(n^{\prime })\oplus \Lambda (n^{\prime \prime }))
\end{equation*}
is identified to a Lagrangian path in $\Lambda (n)$ and 
\begin{equation*}
\bar{\mu}_{\Lambda }(\lambda ^{\prime }\oplus \lambda ^{\prime \prime })=%
\bar{\mu}_{\Lambda ^{\prime }}(\lambda ^{\prime })+\bar{\mu}_{\Lambda
^{\prime \prime }}(\lambda ^{\prime \prime })
\end{equation*}
where $\bar{\mu}_{\Lambda ^{\prime }}$ and $\bar{\mu}_{\Lambda ^{\prime
\prime }}$ are the canonical Lagrangian intersection indices on $\Lambda
(n^{\prime })$ and $\Lambda (n^{\prime \prime })$.
\end{corollary}

\begin{proof}
\smallskip This is trivial by (\ref{foufoun}) and (\ref{muso}).
\end{proof}

The case of the symplectic intersection indices is deduced immediately, 
\textit{mutatis mutandis}:

\begin{corollary}
Let $\sigma ^{\prime }\in C(Sp(n^{\prime }))$ and $\sigma ^{\prime \prime
}\in C(Sp(n^{\prime \prime }))$ be symplectic paths in $Sp(n^{\prime })$ and 
$Sp(n^{\prime \prime })$. Then 
\begin{equation*}
\sigma ^{\prime }\oplus \sigma ^{\prime \prime }\in \mathcal{C}(Sp(n^{\prime
})\oplus Sp(n^{\prime \prime }))
\end{equation*}
is identified with a symplectic path in $Sp$ and 
\begin{equation*}
\bar{\mu}_{Sp}(\sigma ^{\prime }\oplus \sigma ^{\prime \prime })=\bar{\mu}%
_{Sp^{\prime }}(\sigma ^{\prime })+\bar{\mu}_{Sp^{\prime \prime }}(\sigma
^{\prime \prime })
\end{equation*}
where $\bar{\mu}_{Sp^{\prime }}$ and $\bar{\mu}_{Sp^{\prime \prime }}$ are
the canonical symplectic intersection indices on $Sp(n^{\prime })$ and $%
Sp(n^{\prime \prime })$.
\end{corollary}

\subsection{The ``spectral flow'' formula}

Let $(A(t))_{0\leq t\leq 1}$ be a continuous family of real symmetric
matrices of order $n$. We define the \textquotedblleft spectral
flow\textquotedblright\ of this continuous family by the formula 
\begin{equation}
\func{SF}(A(t))_{0\leq t\leq 1}=\func{sign}A(1)-\func{sign}A(0)\text{.}
\label{desf}
\end{equation}%
We have the following result, which has been established in a particular
case by Duistermaat \cite{Duistermaat}.

\begin{proposition}[``Spectral flow formula'']
Let $\lambda _{A}$ be the Lagrangian path defined by 
\begin{equation}
\lambda _{A}(t)=\left\{ x\oplus A(t)x\ ,x\in X\right\}  \label{lama}
\end{equation}
and $\sigma _{A}$ be the symplectic path defined by 
\begin{equation*}
\sigma _{A}(t)=\left( 
\begin{array}{cc}
I & 0 \\ 
A(t) & I%
\end{array}
\right)
\end{equation*}
for $0\leq t\leq 1$. Then 
\begin{equation}
\func{SF}(A(t))_{0\leq t\leq 1}=\bar{\mu}_{\Lambda }(\lambda _{A},X)=\bar{\mu%
}_{Sp}(\sigma _{A},X)  \label{SF}
\end{equation}
where $\func{sign}A(t)$ is the difference between the number of eigenvalues $%
\lambda >0$ and the number of eigenvalues $<0$ of $A(t)$.
\end{proposition}

\begin{proof}
We begin by showing the formula (\ref{SF}). The second one follows
immediately since 
\begin{equation*}
\left( 
\begin{array}{cc}
I & 0 \\ 
A(t) & I%
\end{array}
\right) \left( 
\begin{array}{c}
x \\ 
0%
\end{array}
\right) =\left( 
\begin{array}{c}
x \\ 
A(t)x%
\end{array}
\right) \text{.}
\end{equation*}
Note that $\dim (\lambda _{A}(t)\cap X^{\ast })=n$ for $0\leq t\leq 1$,
hence 
\begin{equation}
\bar{\mu}_{\Lambda }(\lambda _{A},X^{\ast })=0  \label{bonard}
\end{equation}
in view of the axiom (L$_{3}$) which ensures nullity in the strata. In view
of formula (\ref{diff}) in Proposition \ref{chgt} and the antisymmetry of
the signature we have 
\begin{align*}
\bar{\mu}_{\Lambda }(\lambda _{A},X)=\bar{\mu}_{\Lambda }(\lambda _{A},X)-%
\bar{\mu}_{\Lambda }(\lambda _{A},X^{\ast }) \\
=\func{sign}(\lambda _{A}(1),X,X^{\ast })-\func{sign}(\lambda
_{A}(0),X,X^{\ast }) \\
=\func{sign}(X^{\ast },\lambda _{A}(1),X)-\func{sign}(X^{\ast },\lambda
_{A}(0),X)\text{.}
\end{align*}
Since $\func{sign}=\tau $ we have 
\begin{equation*}
\bar{\mu}_{\Lambda }(\lambda _{A},X)=\tau (X^{\ast },\lambda _{A}(0),X)-\tau
(X^{\ast },\lambda _{A}(1),X)
\end{equation*}
and hence 
\begin{equation*}
\bar{\mu}_{\Lambda }(\lambda _{A},X)=\func{sign}A(1)-\func{sign}A(0)
\end{equation*}
in view of property (\ref{taula}) of $\tau $ in Proposition \ref{tautau}.
\end{proof}

This result is evidently rather trivial in the sense that the spectral flow (%
\ref{desf}) depends only on the extreme values $A(1)$ and $A(0)$. The
situation is however far more complicated in the case of infinite
dimensional symplectic spaces and requires elaborated functional analytical
techniques (see \cite{BBF1}).

\section{Comparative study of some other indices}

\subsection{The Robbin-Salamon-McDuff index\label{RSM}}

In \cite{RS1} (see also \cite{Duff}) Robbin and Salamon have constructed,
using differential geometrical methods, a mapping 
\begin{equation*}
\mu _{RS}:\mathcal{C}(\Lambda (n))\times \ell \longrightarrow \tfrac{1}{2}%
\mathbb{Z}
\end{equation*}
which they call ``Maslov index''. They show that their index possesses the
four properties which follow:

\begin{description}
\item[ RS$_{1}$] \textit{If }$\lambda $\textit{\ and }$\lambda ^{\prime }$%
\textit{\ have same endpoints, then }$\mu _{RS}(\lambda ,\ell )=\mu
_{RS}(\lambda ^{\prime },\ell )$ \textit{if and only if }$\lambda \sim
\lambda ^{\prime }$

\item[ RS$_{2}$] $\mu _{RS}(\lambda \ast \lambda ^{\prime },\ell )=\mu
_{RS}(\lambda ,\ell )+\mu _{RS}^{\prime }(\lambda ^{\prime },\ell )$.

\item[ RS$_{3}$] \textit{If} $\func{Im}\lambda \subset \Lambda _{\ell }(n;k) 
$ \textit{then} $\mu _{RS}(\lambda ,\ell )=0$.

\item[ RS$_{4}$] \textit{If }$\gamma $\textit{\ is a loop in }$\Lambda (n)$%
\textit{\ then} $\mu _{RS}(\gamma ,\ell )=m(\gamma )$.

\item[ RS$_{5}$] \textit{If }$\lambda _{A}$\textit{\ is a Lagrangian path (%
\ref{lama}), then}$\ $%
\begin{equation}
\mu _{RS}(\lambda _{A},X)=\tfrac{1}{2}\func{sign}A(1)-\tfrac{1}{2}\func{sign}%
A(0)\text{.}  \label{RSF}
\end{equation}
\end{description}

In fact, there is a simple relation between the indices $\mu _{RS}$ and $%
\bar{\mu}_{\Lambda }$, namely

\begin{proposition}
\label{rikiki}the index $\mu _{RS}$ of Robbin-Salamon is given by 
\begin{equation*}
\mu _{RS}(\lambda ,\ell )=\tfrac{1}{2}\bar{\mu}_{\Lambda }(\lambda ,\ell )
\end{equation*}
where $\bar{\mu}_{\Lambda }$ is the canonical intersection index.
\end{proposition}

\begin{proof}
\smallskip The mapping $2\mu _{RS}$ satisfies the axioms (L$_{1}$--L$_{4}$)
of a Lagrangian intersection index. Therefore, there must exist an index $%
\mu _{\Lambda }$ such that $\mu _{\Lambda }=2\mu _{RS}$. By Theorem \ref{un}
we have thus 
\begin{equation*}
2\mu _{RS}(\lambda ,\ell )=\bar{\mu}_{\Lambda }(\lambda ,\ell )+f(\dim
(\lambda (0)\cap \ell ))-f(\dim (\lambda (1)\cap \ell ))\text{.}
\end{equation*}
for some function $f$. Taking for $\lambda $ the path $\lambda _{A}$ and $%
\ell =X$ we have, in view of the formulae (\ref{SF}) and (\ref{RSF}) : 
\begin{multline*}
\func{sign}A(1)-\func{sign}A(0)=\func{sign}A(1)-\func{sign}A(0)+ \\
f(\dim (\lambda _{A}(0)\cap X))-f(\dim (\lambda _{A}(1)\cap X))\text{.}
\end{multline*}
In view of the equality 
\begin{equation*}
\dim (\lambda _{A}(t)\cap X)=\func{corang}A(t)
\end{equation*}
this is equivalent to 
\begin{equation*}
f(\func{corang}A(0))=f(\func{corang}A(1))
\end{equation*}
for all $A$, Hence $f=0$ since $\func{corang}A(t)$ can take arbitrary values.
\end{proof}

\begin{remark}
Robbin and Salamon use in \cite{RS1} differentiability properties of
Lagrangian paths to construct their index. Furthermore, they invoke deep
functional analytical properties (Kato's theorem concerning the selection of
eigenvalues of a path of symmetric matrices).
\end{remark}

\subsection{The H\"{o}rmander index}

In the frame of his studies of pseudo-differential operators, introduces in 
\cite{OIF} a mapping 
\begin{equation*}
\xi :\Lambda (n)^{4}\ni (\ell _{1},\ell _{2},\ell _{3},\ell
_{4})\longrightarrow \xi (\ell _{1},\ell _{2},\ell _{3},\ell _{4})\in \tfrac{%
1}{2}\mathbb{Z}\text{ .}
\end{equation*}
Robbin and Salamon show that the H\"{o}rmander index is related to their
index $\mu _{RS}$ by the formula 
\begin{equation}
\xi (\ell _{1},\ell _{2},\ell _{3},\ell _{4})=\mu _{RS}(\lambda _{34},\ell
_{2})-\mu _{RS}(\lambda _{34},\ell _{1})  \label{hör}
\end{equation}
where $\lambda _{34}$ is an arbitrary Lagrangian path $\Lambda (n)$ joining $%
\ell _{3}$ to $\ell _{4}$.

\begin{proposition}
The H\"{o}rmander index $\xi $ is given by 
\begin{equation}
\xi (\ell _{1},\ell _{2},\ell _{3},\ell _{4})=\tfrac{1}{2}(\tau (\ell
_{1},\ell _{2},\ell _{3})-\tau (\ell _{1},\ell _{2},\ell _{4}))
\label{ksihö}
\end{equation}%
where $\tau $ is the Demazure--Kashiwara cocycle.
\end{proposition}

\begin{proof}
In view of Proposition \ref{rikiki} and of (\ref{hör}) we have 
\begin{equation*}
\xi (\ell _{1},\ell _{2},\ell _{3},\ell _{4})=\tfrac{1}{2}(\bar{\mu}%
_{\Lambda (n)}(\lambda _{34},\ell _{2})-\bar{\mu}_{\Lambda (n)}(\lambda
_{34},\ell _{1}))\text{.}
\end{equation*}
By the formula (\ref{diff}) in Proposition \ref{chgt} this is 
\begin{equation*}
\xi (\ell _{1},\ell _{2},\ell _{3},\ell _{4})=\tfrac{1}{2}(\func{sign}(\ell
_{4},\ell _{2},\ell _{1})-\func{sign}(\ell _{3},\ell _{2},\ell _{1}))
\end{equation*}
and hence (\ref{ksihö}), since $\func{sign}=\tau $ and the antisymmetry of $%
\tau $.
\end{proof}

\end{document}